\documentclass[graybox]{svmult}

\usepackage{eucal}
\usepackage{amscd}
\usepackage[all]{xy}

\usepackage{mathptmx}       
\usepackage{helvet}         
\usepackage{courier}        
\usepackage{type1cm}        
\usepackage{makeidx}         
\usepackage{graphicx}        
\usepackage{multicol}        
\usepackage[bottom]{footmisc}
\usepackage{cancel}          
\usepackage{dsfont}

\makeindex             

\usepackage{latexsym,amssymb,amsmath} 
\usepackage{mathrsfs} 
\usepackage{color}

\newcommand{\nc}{\newcommand}
\nc{\grad}[1]{^{({#1})}} \nc{\fil}[1]{_{#1}}
\nc{\surj}{\to\hskip -3mm \to}
\nc{\mop}[1]{\mathop{\hbox {\rm #1} }\nolimits}
\nc{\smop}[1]{\mathop{\hbox {\eightrm #1} }\nolimits}
\nc{\mopl}[1]{\mathop{\hbox {\rm #1} }\limits}
\nc{\un}{\hbox{\bf 1}}

\nc{\cala}{{\mathcal A}}
\nc{\calb}{{\mathcal B}}
\nc{\calc}{{\mathcal C}}
\nc{\cald}{{\mathcal D}}
\nc{\cale}{{\mathcal E}}
\nc{\calf}{{\mathcal F}}
\nc{\calg}{{\mathcal G}}
\nc{\calh}{{\mathcal H}}
\nc{\cali}{{\mathcal I}}
\nc{\calj}{{\mathcal J}}
\nc{\call}{{\mathcal L}}
\nc{\calm}{{\mathcal M}}
\nc{\caln}{{\mathcal N}}
\nc{\calo}{{\mathcal O}}
\nc{\calp}{{\mathcal P}}
\nc{\calr}{{\mathcal R}}
\nc{\calt}{{\mathcal T}}
\nc{\calu}{{\mathcal U}}
\nc{\calv}{{\mathcal V}}
\nc{\calw}{{\mathcal W}}
\nc{\calx}{{\mathcal X}}

\allowdisplaybreaks[4]

\definecolor{Light}{gray}{0.85}

\def\abs#1{\left\vert #1 \right\vert}
\def\allpoly{\mbox{$\re\langle X \rangle$}}

\def\allpolyx0degn{\mbox{$P_n$}}

\def\allseries{\mbox{$\re\langle\langle X \rangle\rangle$}}

\def\allseriesdeltamLC{\mbox{$\re^m_{LC}\langle\langle X_\delta \rangle\rangle$}}
\def\allseriesdeltam{\mbox{$\re^m\langle\langle X_\delta \rangle\rangle$}}
\def\allseriesell{\mbox{$\re^{\ell} \langle\langle X \rangle\rangle$}}

\def\allseriesLC{\mbox{$\re_{LC}\langle\langle X \rangle\rangle$}}
\def\allseriesm{\mbox{$\re^m\langle\langle X \rangle\rangle$}}

\def\allseriesmLC{\mbox{$\re^{m}_{LC}\langle\langle X \rangle\rangle$}}

\def\allseriesellLC{\mbox{$\re^{\ell}_{LC}\langle\langle X \rangle\rangle$}}
\def\allseriesellGC{\mbox{$\re^{\ell}_{GC}\langle\langle X \rangle\rangle$}}

\def\bfem#1{{\bf \em #1}} 

\newcommand{\comment}[1]{} 

\def\Endallseries{{\rm End}(\allseries)}

\def\eqref#1{(\ref{#1})} 

\def\expup{{\rm e}}

\def\Fliessdelta{\mathscr{F}_{\delta}}

\def\id{{\rm id}}

\def\Lpm{L_{\mathfrak{p}}^m}

\def\Lpme{L^m_{\mathfrak{p},e}}

\def\liepoly{{\cal L}(X)}

\def\modcomp{\:\tilde{\circ}\,} 

\def\nat{{\mathbb N}} 
\def\norm#1{\left\Vert#1\right\Vert}

\def\openbull{\framebox[0.08in][c]{$\;$}} 
\def\ord{{\rm ord}}

\newcommand{\R}{\mathbb R}

\def\re{{\mathbb R}} 

\def\sameau{\rule[0.017in]{0.2in}{0.012in}}

\def\shuffle{{\scriptscriptstyle \;\sqcup \hspace*{-0.05cm}\sqcup\;}}

\def\supp{{\rm supp}}

\def\begals{\[\begin{aligned}}
\def\endals{\end{aligned}\]}
\def\begce{\begin{center}}
\def\endce{\end{center}}
\def\begar{\begin{array}}
\def\endar{\end{array}}
\def\begeq{\begin{equation}}
\def\endeq{\end{equation}}
\def\begdi{\begin{displaymath}}
\def\enddi{\end{displaymath}}
\def\begdis{\begin{eqnarray*}}
\def\enddis{\end{eqnarray*}}
\def\begeqa{\begin{eqnarray}}
\def\endeqa{\end{eqnarray}}
\def\begdes{\begin{description}}
\def\enddes{\end{description}}
\def\begit{\begin{itemize}}
\def\endit{\end{itemize}}
\def\begen{\begin{enumerate}}
\def\enden{\end{enumerate}}
\def\beglar{\left[\begin{array}}
\def\endrar{\end{array}\right]}
\def\begle{\begin{lemma}}
\def\endle{\end{lemma}}
\def\begde{\begin{definition}}
\def\endde{\end{definition}}
\def\begth{\begin{theorem}}
\def\endth{\end{theorem}}
\def\begco{\begin{corollary}}
\def\endco{\end{corollary}}
\def\begprop{\begin{proposition}}
\def\endprop{\end{proposition}}
\def\begex{\begin{example}}
\def\endex{\hfill\openbull \end{example} \vspace*{0.15in}}
\def\begexer{\begin{exercise}}
\def\endexer{\end{exercise}}

\def\begres{\noindent{\bf Remarks}:\begin{enumerate}}
\def\endres{\end{enumerate} \par}

\def\begtab{\begin{tabular}}
\def\endtab{\end{tabular}}
\def\rref#1{(\ref{#1})}
\def\lbracedef{\left\{\begin{array}{@{\hspace*{2pt}}c@{\hspace*{3pt}}c@{\hspace*{3pt}}l}} 

\def\allpolybarX{\mbox{$\re \langle \bar{X} \rangle$}}
\def\allseriesbarX{\mbox{$\re \langle\langle \bar{X} \rangle\rangle$}}

\def\allseriesXc{\mbox{$\re \langle\langle X_{c} \rangle\rangle$}}
\def\allseriesXd{\mbox{$\re \langle\langle X_{d} \rangle\rangle$}}
\def\allseriesXe{\mbox{$\re \langle\langle X_{e} \rangle\rangle$}}

\def\allseriesXRat{\mbox{$\re_{rat} \langle\langle X \rangle\rangle$}}

\def\allseriesXeRat{\mbox{$\re_{rat} \langle\langle X_{e} \rangle\rangle$}}

\begin{document}

\title*{Combinatorial Hopf algebras for interconnected nonlinear input-output
systems with a view towards discretization}

\titlerunning{Combinatorial Hopf algebra for interconnected nonlinear systems}

\author{Luis A.~Duffaut Espinosa, Kurusch Ebrahimi-Fard, and W.~Steven Gray}
\institute{Luis A.~Duffaut Espinosa
\at Department of Electrical and Biomedical Engineering,
University of Vermont, Burlington, Vermont 05405, USA,
\email{l.duffautespinosa@gmail.com}
\url{}
\and Kurusch Ebrahimi-Fard
\at Department of Mathematical Sciences,
Norwegian University of Science and Technology -- NTNU,
NO-7491 Trondheim, Norway,
\email{kurusch.ebrahimi-fard@ntnu.no},
\url{https://folk.ntnu.no/kurusche/}
\and W.~Steven Gray
\at Department of Electrical and Computer Engineering,
Old Dominion University, Norfolk, Virginia 23529, USA,
\email{sgray@odu.edu}
\url{}}

%
%

\maketitle

\abstract*{A detailed expose of the Hopf algebra approach to interconnected input-output systems in nonlinear control theory is presented. The focus is on input-output systems that can been represented in terms of Chen--Fliess functional expansions or Fliess operators. This provides a starting point for a discrete-time version of this theory. In particular, the notion
of a discrete-time Fliess operator is given and a class of parallel interconnections is described.}

\abstract{A detailed expose of the Hopf algebra approach to interconnected input-output systems in nonlinear control theory is presented. The focus is on input-output systems that can been represented in terms of Chen--Fliess functional expansions or Fliess operators. This provides a starting point for a discrete-time version of this theory. In particular, the notion of a discrete-time Fliess operator is given and a class of parallel interconnections is described in terms of the quasi-shuffle algebra.}


\section{Introduction}
\label{sec:introduction}

A central problem in control theory is understanding how dynamical systems
behave when they are interconnected. In a typical design problem, one is given a
fixed system representing the {\em plant}, say a robotic manipulator or a spacecraft.
The objective is to find a second
system, usually called the {\em controller},
which when interconnected with the first will make the output track a
pre-specified trajectory. When the component systems are nonlinear, these
problems are difficult to address by purely analytical means. The prevailing
methodologies are geometric in nature and based largely on state variable
analysis \cite{Isidori_95,Khalil_14,Nijmeijer-vanderSchaft_90,Slotine-Li_91}.
A complementary approach, however, has begun to
emerge where the input-output map of each component system is represented in
terms of a Chen--Fliess functional expansion or Fliess operator. In this setting,
concepts from combinatorics and algebra are employed to produce an explicit
description of the interconnected system. The genesis of this method can be found in the work
of Fliess \cite{Fliess_81,Fliess_83} and Ferfera \cite{Ferfera_79,Ferfera_80},
who described key elements of the underlying algebras in terms of noncommutative
formal power series. In particular, they identified the central role played by
the shuffle algebra in this theory. The method was further developed by Gray et
al.~\cite{Gray-Li_05,Gray-Thitsa_IJC12,Gray-Wang_SCL02,Thitsa-Gray_12} and Wang
\cite{Wang_90} who addressed basic analytical questions such as what inputs
guarantee convergence of the component series, when are the interconnections
well defined, and what is the nature of the output functions? A fundamental open
problem on the algebraic side until 2011 was how to explicitly compute the
generating series of two feedback interconnected Fliess operators. Largely inspired by
interactions at the 2010 Trimester in Combinatorics and Control
(COCO2010) in Madrid with researchers in the area of quantum field theory (see, for
example, \cite{Figueroa-Gracia-Bondia_05}), the problem was eventually solved by
identifying a new combinatorial Hopf algebra underlying the calculation.
The evolution of this structure took several distinct steps.
The
single-input, single-output (SISO) case was first addressed by Gray and Duffaut
Espinosa in \cite{Gray-Duffaut_Espinosa_SCL11} via a certain graded Hopf algebra
of combinatorial type. Foissy then introduced a new grading in \cite{Foissy_CiA15}
which rendered a connected version of this combinatorial Hopf algebra. This
naturally provided a fully recursive formula for the antipode, which is central
to the feedback calculation \cite{Gray-et-al_MTNS14}. The multivariable, i.e.,
multi-input, multi-output (MIMO) case, was then treated in
\cite{Gray-etal_SCL14}. Next, a full combinatorial treatment, including a
Zimmermann type forest formula for the antipode \cite{Anshelevich-etal_05}, was presented in
\cite{Duffaut-etal_submit2014}. This last result, based on an equivalent
combinatorial Hopf algebra of decorated rooted circle trees, greatly reduces the
number of computations involved by eliminating the inter-term cancelations that
are intrinsic in the usual antipode calculation. Practical problems would be
largely intractable without this innovation.  The final and most recent development method is a description of this Hopf algebra based entirely on (co)derivation(-type) maps applied to the (co)product. This method was first observed to be implicit in the work of Devlin on the classical Poincar\'{e} center problem \cite{Devlin_90,Ebrahimi-Fard-Gray_IMRN17}. In a state space setting, it can be related to computing iterated Lie derivatives to determine series coefficients \cite{Isidori_95}. Control applications ranging from guidance and chemical engineering to systems biology can be found in \cite{Duffaut_Espinsa-Gray_ACC17,Gray-Duffaut_Espinosa_IRMA15,Gray-etal_SCL14,Gray-etal_CISS15,Gray-etal_CDC15,Gray-et-al_AUT14}.

\smallskip

This article has two general goals. First, an introduction to the method of
combinatorial Hopf algebra in the context of feedback control theory is given
in its most complete and up-to-date form. The idea is to integrate all of the advances
described above into a single uniform treatment. This results in a new and distinct
presentation of these ideas. In particular, the full solution to the problem of computing the
generating series for a multivariable continuous-time dynamic output feedback
system will be described. The origins of the forest formula
related to this calculation will also be outlined. Then the focus is shifted to
the second objective, which is largely an open problem in this field, namely how to
recast this theory in a discrete-time setting. One approach suggested by Fliess in
\cite{Fliess_80} is to describe the generating series of a discrete-time
input-output map in terms of a complete tensor algebra defined on an algebra of
polynomials over a noncommutative alphabet. While this is a very general approach and can be
related to the notion of a Volterra series as described by Sontag in \cite{Sontag_79},
to date it has not lead to any particularly useful algebraic structures in the context
of system interconnections.
Another approach developed by the authors in the context of numerical approximation
is to define the notion of a discrete-time Fliess operator in terms of a series of
iterated sums over a noncommutative alphabet \cite{Gray-etal_ACC16,Gray-etal_NM}.
While not the most general set up, it has been shown to be related to a class
of state affine rational input discrete-time systems in the case where the generating series
is rational. Furthermore, this class of so called
rational discrete-time Fliess operators is guaranteed to always converge.
So this will be the approach taken here.  As the main interest in this paper is interconnection theory, the analysis
begins with the simplest type of interconnections, the parallel sum and
parallel product connections. The former is completely trivial, but the later
induces the quasi-shuffle algebra of Hoffman (see \cite{Hoffman_00}) on the vector space of
generating series. Of particular interest is whether rationality is preserved under the
quasi-shuffle product. It is well known to be the case for the shuffle product \cite{Fliess_81}.

\smallskip

The paper is organized as follows. In Section \ref{sec:preliminaries}, some preliminaries
on Fliess operators and graded connected Hopf algebras are given to set the notation and
provide some background. The subsequent section is devoted to describing the combinatorial
algebras that are naturally induced by the interconnection of Fliess operators.
In Section~\ref{sec:cha-crt}, the Hopf algebra of coordinate functions for the output feedback
Hopf algebra is described in detail. The final section addresses elements of the discrete-time version
of this theory.


\section{Preliminaries}
\label{sec:preliminaries}

A finite nonempty set of noncommuting symbols $X=\{ x_0,x_1,$ $\ldots,x_m\}$ is
called an {\em alphabet}. Each element of $X$ is called a {\em letter}, and any
finite sequence of letters from $X$, $\eta=x_{i_1}\cdots x_{i_k}$, is called a
{\em word} over $X$. The {\em length} of the word $\eta$, denoted $\abs{\eta}$, is given
by the number of letters it contains. The set of all words with length $k$ is
denoted by $X^k$. The set of all words including the empty word, $\emptyset$, is
designated by $X^\ast$, while $X^+:=X^\ast-\{\emptyset\}$. The set $X^\ast$
forms a monoid under catenation. The set $\eta X^\ast$ is comprised of all words
with the prefix $\eta\in X^\ast$. For any fixed integer $\ell\geq 1$, a mapping
$c:X^\ast\rightarrow \re^\ell$ is called a {\em formal power series}. The value
of $c$ at $\eta\in X^\ast$ is written as $(c,\eta)\in\re^{\ell}$ and called the
{\em coefficient} of the word $\eta$ in $c$. Typically, $c$ is represented as
the formal sum $c=\sum_{\eta\in X^\ast}(c,\eta)\eta.$ If the {\em constant term}
$(c,\emptyset)=0$ then $c$ is said to be {\em proper}. The {\em support} of $c$,
$\supp(c)$, is the set of all words having nonzero coefficients in $c$. The {\em
order} of $c$, $\ord(c)$, is the length of the shortest word in its support
($\ord(0):=\infty$).\footnote{For notational convenience,
$p=(p,\emptyset)\emptyset\in\allpoly$ is often abbreviated as
$p=(p,\emptyset)$.} The {$\re$-vector space}
of all formal power series over $X^*$ with coefficients in $\mathbb{R}^\ell$ is
denoted by $\allseriesell$.\footnote{{The superscript $\ell$ will be dropped when $\ell=1$.}}
It forms a unital associative $\re$-algebra under
the catenation product. {Specifically, for $c,d\in\allseriesell$ the catenation
product is given by $cd=\sum_{\eta\in X^\ast} (cd,\eta)\,\eta$, where
\begdi
(cd,\eta)=\sum_{\eta=\xi\nu} (c,\xi)(d,\nu), \;\;\forall\eta\in X^{\ast},
\enddi
and the product on $\re^\ell$ is defined componentwise.
The unit in this case is $\mathbf{1}:=1\emptyset$.}
$\allseriesell$ is also a unital, commutative and associative $\re$-algebra
under the shuffle product, denoted here by the shuffle symbol $\shuffle$. The
shuffle product of two words is defined inductively by
\begin{equation}
\label{shuffleproduct}
	(x_i\eta) \shuffle (x_j\xi)=x_i(\eta\shuffle(x_j\xi))+x_j((x_i\eta)\shuffle \xi)
\end{equation}
with $\eta\shuffle\emptyset=\emptyset\shuffle\eta=\eta$ given any words
$\eta,\xi\in X^\ast$ and letters $x_i,x_j\in X$ \cite{Fliess_81,Ree_58,Reutenauer_93}.
For instance, $x_i \shuffle x_j = x_ix_j + x_jx_i$ and
\begin{equation*}
	x_{i_1} x_{i_2}\shuffle x_{i_3} x_{i_4} = x_{i_1} x_{i_2}x_{i_3} x_{i_4}
	+ x_{i_3} x_{i_4}x_{i_1} x_{i_2} + x_{i_1} x_{i_3}(x_{i_2}\shuffle
x_{i_4})
	+ x_{i_3} x_{i_1}(x_{i_2}\shuffle x_{i_4}).
\end{equation*}
{The definition of the shuffle product
is extended linearly to any two series $c,d\in \allseriesell$ by letting
\begeq \label{eq:shuffle}
c \shuffle d = \sum_{\eta,\xi\in X^{\ast}} (c,\eta)(d,\xi)\,\eta\shuffle\xi,
\endeq
where again the product on $\re^\ell$ is defined componentwise.
For a fixed word $\nu\in X^{\ast}$, the coefficient
\begdi
(\eta\shuffle \xi,\nu)=0\;\; {\rm if}\;\; |\eta|+|\xi|\neq |\nu|.
\enddi
Hence, the infinite sum in (\ref{eq:shuffle}) is always well defined since the family
of polynomials
$\{\eta\shuffle\xi\}_{\eta,\xi\in X^{\ast}}$ is locally finite. The unit for this product is
$\mathbf{1}$.
}


Some standard concepts regarding rational formal power series, which are used in Section \ref{sec:towards-discretization}, are provided next
{\cite{Berstel-Reutenauer_88}}.
A series $c\in\allseries$ is called {\em
invertible} if there exists a series $c^{-1}\in\allseries$ such that
$cc^{-1}=c^{-1}c=\mathbf{1}$. 
In the
event that $c$ is not proper, it is always possible to write
\begdi
	c=(c,\emptyset)(\mathbf{1}-c^{\prime}),
\enddi
where $(c,\emptyset)$ is nonzero, and $c^{\prime}\in\allseries$ is proper.
It then follows that
\begeq \label{eq:Kleene_inverse}
	c^{-1}=\frac{1}{(c,\emptyset)}(\mathbf{1}-c^{\prime})^{-1}
	         =\frac{1}{(c,\emptyset)} (c^\prime)^{\ast},
\endeq
where
\begdi
	(c^{\prime})^{\ast}:=\sum_{i=0}^{\infty} (c^{\prime})^i.
\enddi
In fact, $c$ is invertible if and {\em only if} $c$ is not proper.
Now let $S$ be a subalgebra of the $\re$-algebra $\allseries$
with the catenation product. $S$ is said to be {\em rationally closed} when every invertible
$c\in S$ has $c^{-1}\in S$ (or equivalently, every proper $c^{\prime}\in S$ has $(c^{\prime})^{\ast}\in
S$). The {\em rational closure} of any
subset $E\subset\allseries$ is the smallest rationally closed
subalgebra of $\allseries$ containing $E$.

\begde  \cite{Berstel-Reutenauer_88}  \label{def:rationality}
A series $c\in\allseries$ is \bfem{rational} if it belongs to the
rational closure of $\allpoly$. The subset of all rational series in
$\allseries$  is denoted by $\allseriesXRat$.
\endde

{From Definition \ref{def:rationality} it is clear that every series
in $\allseriesXRat$ is generated by a finite number of rational operations (scalar multiplication, addition, catenation, and inversion) applied to a finite
set of polynomials over $X$. In the case where the alphabet $X$ has an infinite number of letters,
such a series can only involve a finite subset of $X$. Therefore, it is rational in exactly the sense described above when restricted to this sub-alphabet.
This will be the notion of rationality employed in this manuscript whenever $X$ is infinite. But the reader is cautioned that
other notions of rationality for infinite alphabets appear in the literature, see, for example, \cite{Pittou-Rahonis_2017}.}

It turns out that an entirely different characterization of a rational series
is possible using
{a monoid structure on the set of
$n\times n$ matrices over $\re$, denoted $\re^{n\times n}$, where the product is
conventional matrix multiplication and the unit is the $n\times n$ identity matrix $I$.}

\begde  \cite{Berstel-Reutenauer_88} \label{def:linear-representation}
A \bfem {linear} \bfem{representation}
of a series $c\in\allseries$ is any triple
$(\mu,\gamma,\lambda)$, where
\begdi
	\mu: X^{\ast}\rightarrow \re^{n\times n}
\enddi
is a monoid morphism, and $\gamma,\lambda^T\in\re^{n\times 1}$ are such that
\begeq \label{eq:coefs-from-representation}
	(c,\eta)=\lambda\mu(\eta)\gamma, \quad \forall\eta\in X^{\ast}.
\endeq
The integer $n$ is the dimension of the representation.
\endde

\begde {\cite{Berstel-Reutenauer_88}}
A series $c\in\allseries$ is called \bfem{recognizable} if
it has a linear representation.
\endde

\begth {{\rm\cite{Schutzenberger_61}}} \label{th:Schutzenberger}
A formal power series is rational if and only if it is recognizable.
\endth

A third characterization of rationality is given by the notion of {\em stability}.
{Define for any letter $x_i\in X$ and word $\eta=x_j \eta' \in X^\ast$ the left-shift operator
\begeq \label{eq:left_shift}
x_i^{-1}(\eta) = \delta_{ij} \eta',
\endeq
where $\delta_{ij}$ is the standard Kronecker delta}.
Higher order shifts are defined inductively via
$(x_i\xi)^{-1}(\cdot)=\xi^{-1}x_i^{-1}(\cdot)$, where $\xi\in X^\ast$. The
left-shift operator is assumed to act linearly on $\allseries$.

\begde  \cite{Berstel-Reutenauer_88}
A subset $V\subset\allseries$ is called \bfem{stable} when $\xi^{-1}(c)\in V$
for all $c\in V$ and $\xi\in X^{\ast}$.
\endde

\begth  {\rm \cite{Berstel-Reutenauer_88}} \label{th:stable-rational-subspace}
A series $c\in\allseries$ is rational if and only if
there exists a stable finite dimensional $\re$-vector subspace of
$\allseries$ containing $c$.
\endth


\subsection{{Chen--Fliess series and} Fliess operators}
\label{subsec:Fliess-operators}

{One can associate with any formal power series $c\in\allseriesell$
a functional series, $F_c$, known a {\em Chen--Fliess series}. }
Let $\mathfrak{p}\ge 1$ and $t_0 < t_1$ be given. For a Lebesgue measurable
function $u: [t_0,t_1] \rightarrow\re^m$, define
$\norm{u}_{\mathfrak{p}}=\max\{\norm{u_i}_{\mathfrak{p}}: \ 1\le
i\le m\}$, where $\norm{u_i}_{\mathfrak{p}}$ is the usual
$L_{\mathfrak{p}}$-norm for a measurable real-valued function,
$u_i$, defined on $[t_0,t_1]$.  Let $L^m_{\mathfrak{p}}[t_0,t_1]$
denote the set of all measurable functions defined on $[t_0,t_1]$
having a finite $\norm{\cdot}_{\mathfrak{p}}$ norm and
$B_{\mathfrak{p}}^m(R)[t_0,t_1]:=\{u\in
L_{\mathfrak{p}}^m[t_0,t_1]:\norm{u}_{\mathfrak{p}}\leq R\}$. Assume
$C[t_0,t_1]$ is the subset of continuous functions in $L_{1}^m[t_0,t_1]$.
Define inductively for each word $\eta\in X^{\ast}$ the map $E_\eta:
L_1^m[t_0, t_1]\rightarrow C[t_0, t_1]$ by setting $E_\emptyset[u]=1$ and
letting
\begin{align*}
	E_{x_{i_1}\bar{\eta}}[u](t,t_0) &:=
	\int_{t_0}^tu_{i_1}(\tau_1)E_{\bar{\eta}}[u](\tau_1,t_0)\,d\tau_1\\
	&= \int_{\Delta^n_{[t,t_0]}} u_{i_1}(\tau_1)u_{i_2}(\tau_2) \cdots
	u_{i_n}(\tau_n) d\tau_n \cdots d\tau_{2} d\tau_1,
\end{align*}
where $\Delta^n_{[t,t_0]}:=\{(\tau_1,\ldots, \tau_n),\  t\ge \tau_1 \ge \cdots \ge
\tau_n \ge t_0\},$ $\eta = x_{i_1} \cdots x_{i_n} = x_{i_1} \bar{\eta} \in
X^{\ast}$, and $u_0:=1$. For instance, the words $x_i$ and $x_{i_1}x_{i_2}$
correspond to the integrals
$$
	E_{x_i}[u](t,t_0) =
	\int_{t_0}^tu_{i}(\tau)d\tau,
	\quad
	E_{x_{i_1}x_{i_2}}[u](t,t_0) =
	\int_{t_0}^tu_{i_1}(\tau_1)\int_{t_0}^{\tau_1}u_{i_2}(\tau_2)d\tau_2d\tau_1.
$$
The {{\em Chen--Fliess series} corresponding to $c\in\allseriesell$ is defined} to be
\begeq
	F_c[u](t) = \sum_{\eta\in X^{\ast}} (c,\eta)\,E_\eta[u](t,t_0)
\label{eq:Fliess-operator-defined}.
\endeq
In the event that
there exist real numbers $K_c,M_c>0$ such that
\begin{equation}
	\abs{(c,\eta)}\le K_c M_c^{|\eta|}|\eta|!,\quad \forall\eta\in X^{\ast},
\label{eq:local-convergence-growth-bound}
\end{equation}
then {$F_c$ constitutes a well defined causal operator} from $B_{\mathfrak p}^m(R)[t_0,$
$t_0+T]$ into $B_{\mathfrak q}^{\ell}(S)[t_0, \, t_0+T]$ for some $S>0$ provided
$\bar{R}:=\max\{R,T\}<1/M_c(m+1)$,
and the numbers $\mathfrak{p},\mathfrak{q}\in[1,\infty]$ are conjugate
exponents, i.e., $1/\mathfrak{p}+1/\mathfrak{q}=1$ \cite{Gray-Wang_SCL02}.
(Here, $\abs{z}:=\max_i \abs{z_i}$ when $z\in\re^\ell$.) In this case,
{$F_c$ is called a {\em Fliess operator}}
and said to be {\em locally convergent} (LC). The set of all
series satisfying \rref{eq:local-convergence-growth-bound} is denoted by
$\allseriesellLC$. When $c$ satisfies the more stringent growth condition
\begeq
\abs{(c,\eta)}\le K_c M_c^{|\eta|},\quad \forall\eta\in X^{\ast},
\label{eq:global-convergence-growth-bound}
\endeq
the series \rref{eq:Fliess-operator-defined}
defines a {Fliess operator} from the extended space
$\Lpme (t_0)$ into $C[t_0, \infty)$,  where
\begin{align*}
\Lpme(t_0):=
\{u:[t_0,\infty)\rightarrow \re^m:u_{[t_0,t_1]}\in \Lpm[t_0,t_1],\quad \forall t_1
\in (t_0,\infty)\},
\end{align*}
and $u_{[t_0,t_1]}$ denotes the restriction of $u$ to the interval $[t_0,t_1]$
\cite{Gray-Wang_SCL02}. In this case, the operator is said to be {\em globally
convergent} (GC), and the set of all series satisfying
\rref{eq:global-convergence-growth-bound} is designated by
$\allseriesellGC$.

{
Most of the work regarding nonlinear input-output systems in
control theory prior to the work of Fliess was
based on Volterra series, see, for example, \cite{Brockett_76,Lesiak-Krener_78,Rugh_81}.
In many ways this earlier work set the stage for the introduction of the
noncommutative algebraic framework championed by Fliess.
}
As Fliess operators are series of weighted iterated integrals of control
functions, they are also related to
the work of K.~T.~Chen, who revealed that iterated integrals come with a natural
algebraic structure \cite{Chen_54,Chen_57,Chen_67b}. Indeed, products of iterated integrals can again be written
as linear combinations of iterated integrals. This is implied by the classical
integration by parts rule for indefinite Riemann integrals, which yields for
instance that
\begin{align*}
	E_{x_{i_1}}[u](t,t_0) E_{x_{i_2}}[u](t,t_0) =
	E_{x_{i_1}x_{i_2}}[u](t,t_0) + E_{x_{i_2}x_{i_1}}[u](t,t_0).
\end{align*}	
Linearity allows one to relate this to the shuffle product \eqref{shuffleproduct}
$$
	E_{x_{i_1}}[u](t,t_0) E_{x_{i_2}}[u](t,t_0) =
	F_{x_{i_1}x_{i_2}+x_{i_2}x_{i_1}}[u](t)
	= F_{x_{i_1}\shuffle x_{i_2}}[u](t).
$$
This generalizes naturally to the shuffle product for iterated integrals with
respects to words $\eta,\nu \in X^*$
\begeq
\label{ShuffleProduct1}
	E_{\eta}[u](t,t_0) E_{\nu}[u](t,t_0) = F_{\eta \shuffle \nu}[u](t),
\endeq
which in turn implies for Fliess operators corresponding to $c,d \in
\allseriesell$ that
\begeq
\label{ShuffleProduct2}
	F_c[u](t) F_d[u](t) = F_{c \shuffle d}[u](t).
\endeq

{
A Fliess operator $F_c$ defined on $B_{\mathfrak p}^m(R)[t_0,t_0+T]$
is said to be {\em realizable} when there exists
a state space model consisting of $n$ ordinary differential
equations and $\ell$ output functions
\begin{subequations}
\begin{align}
\dot{z}(t)&= g_0(z(t))+\sum_{i=1}^m g_i(z(t))\,u_i(t),\;\;z(t_0)=z_0  \label{eq:state}\\
y_j(t)&=h_j(z(t)),\;\;j=1,2,\ldots,\ell, \label{eq:output}
\end{align}
\end{subequations}
where each $g_i$ is an analytic vector field expressed in local
coordinates on some neighborhood ${\cal W}$ of $z_0$,
and each output function $h_j$ is an analytic function on ${\cal W}$ such that
(\ref{eq:state}) has a well defined solution $z(t)$, $t\in[t_0,t_0+T]$ for
any given input $u\in B_{\mathfrak
p}^m(R)[t_0,t_0+T]$, and $y_j(t)=F_{c_j}[u](t)=h_j(z(t))$, $t\in[t_0,t_0+T]$, $j=1,2,\ldots,\ell$.
It can be shown that for any word $\eta=x_{i_k}\cdots x_{i_1}\in X^\ast$
\begeq
(c_j,\eta)=L_{g_{\eta}}h_j(z_0):=L_{g_{i_1}}\cdots L_{g_{i_k}}h_j(z_0), \label{eq:c-equals-Lgh}
\endeq
where $L_{g_i}h_j$ is the {\em Lie derivative} of $h_j$ with respect to $g_i$.
For any $c\in\allseriesell$, the $\re$-linear mapping
${\cal H}_c:\allpoly\rightarrow \allseriesell$
uniquely specified by
$
({\cal H}_c(\eta),\xi)=(c,\xi\eta),\;\; \xi,\eta \in X^{\ast}
$
is called the {\em Hankel mapping} of $c$. The series $c$ is said to
have finite {\em Lie rank} $\rho_L(c)$ when the range of ${\cal H}_c$
restricted to the $\re$-vector space of Lie polynomials over $X$,
i.e., the free Lie algebra $\liepoly\subset \allpoly$, has dimension $\rho_L(c)$.
It is well known that $F_c$ is realizable if and
only if $c\in\allseriesellLC$ has finite Lie rank
\cite{Fliess_81,Fliess_83,Isidori_95,Jakubczyk_80,Jakubczyk_86,Jakubczyk_86a,Sussmann_77,Sussmann_90}. In which case, all minimal realization
have dimension $\rho_L(c)$ and are unique up to a diffeomorphism.
In the event that ${\cal H}_c$ has
finite rank on the entire vector space $\allpoly$, usually referred to as the {\em Hankel rank} $\rho_H(c)$ of $c$,
and $c\in\allseriesellGC$ then
$F_c$ has a minimal bilinear state space realization
\begin{align*}
\dot{z}(t) &=A_0z(t)+\sum\limits_{i=1}^m A_iz(t)u_i(t),\;\;z(t_0)=z_0 \\
 y_j(t) &=C_j z(t),\;\;j=1,2,\ldots,\ell
\end{align*}
of dimension $\rho_H(c)\geq \rho_L(c)$,
where $A_j$ and $C_j$ are real matrices of appropriate dimensions \cite{Fliess_81,Fliess_83,Isidori_95}. Here
the state $z(t)$ is well defined on any interval $[t_0,t_0+T]$, $T>0$, when $u\in L_{1,e}^m(t_0)$,
and the operator $F_c$ always converges globally.
In addition, \rref{eq:c-equals-Lgh} simplifies to
\begeq
(c_j,x_{i_k}\cdots x_{i_1})=C_jA_{i_k}\cdots A_{i_1}z_0. \label{eq:c-equals-Lgh-bilinear}
\endeq
In light of \rref{eq:coefs-from-representation}, it is not hard to see that $c$ is recognizable in SISO case if and only if
$F_c$ has a bilinear realization with $A_i=\mu(x_i)$ for $i=0,1$ and $z_0=\gamma$, and $C=\lambda$.
}


\subsection{Graded connected Hopf algebras}
\label{subsec:Hopf-algebras}

All algebraic structures here are considered over the base field $\mathbb{K}$ of characteristic zero, for
instance, $\mathbb{C}$ or  $\mathbb{R}$. Multiplication in $\mathbb{K}$ is
denoted by $m_\mathbb{K}: \mathbb{K} \otimes \mathbb{K} \to \mathbb{K}$.


\subsubsection{Algebra}
\label{ssubsec:algebras}

A {\it $\mathbb{K}$-algebra} is denoted by the triple $(A,m_A,\eta_A)$ where $A$
is a $\mathbb{K}$-vector space carrying an associative product $m_A: A \otimes A
\to A$, i.e., $m_A \circ (m_A \otimes \mathrm{id}_A) = m_A \circ (\mathrm{id}_A
\otimes m_A): A \otimes A \otimes A \to A$, and a unit map $\eta_A: \mathbb{K}
\to A$. The algebra unit corresponding to the latter is denote by $1_A$. A {\it
$\mathbb{K}$-subalgebra} of the $\mathbb{K}$-algebra $A$ is a
$\mathbb{K}$-vector subspace $B \subseteq A$ such that $m_A(b \otimes b') \in B$
for all $b,b' \in B$. A $\mathbb{K}$-subalgebra $I \subseteq A$ is called a  {\it
(right-) left-ideal} if for any elements $i \in I$ and $a \in A$ the product
($m_A( i \otimes a)$) $m_A(a \otimes i)$ is in $I$. An  {\it ideal} $I \subseteq
A$ is both a left- and right-ideal.

In order to motivate the concept of a $\mathbb{K}$-coalgebra, the definition of a
$\mathbb{K}$-algebra $A$ is rephrased in terms of commutative diagrams.
Associativity of the $\mathbb{K}$-vector space morphism $m_A: A \otimes A \to A$
translates into commutativity of the diagram
\begin{equation}
    \xymatrix{ A \otimes A \otimes A \ar[rr]^{m_A \otimes \mathrm{id}_A}
       \ar[d]_{\mathrm{id}_A \otimes m_A} && A \otimes A \ar[d]^{m_A}\\
        A \otimes A \ar[rr]^{m_A} && A }
    \label{assoc2}
\end{equation}
The $\mathbb{K}$-algebra $A$ is unital if the $\mathbb{K}$-vector space map
$\eta_A: \mathbb{K}
\to A$ satisfies the commutative diagram
\begin{equation}
 \xymatrix{ \mathbb{K} \otimes A \ar[r]^{\eta_A \otimes \mathrm{id}_A}
 \ar[rd]_{\alpha_l} & A \otimes A \ar[d]^{m_A} & A \otimes \mathbb{K}
 \ar[l]_{\mathrm{id}_A \otimes \eta_A}
 \ar[ld]^{\alpha_r} \\
  & A & }
    \label{unit}
\end{equation}
Here $\alpha_l$ and $\alpha_r$ are the isomorphisms sending $k \otimes a$
respectively $a \otimes k$ to $ka$ for $k \in  \mathbb{K}$, $a \in A$. Let
$\tau:= \tau_{A,A}: A \otimes A \to A \otimes A$ be the flip map, $\tau_{A,A}(x
\otimes y):= y \otimes x$. The $\mathbb{K}$-algebra $A$ is {\it commutative} if
the next diagram commutes
\begin{equation}
 \xymatrix{
     A \otimes A \ar[d]_{m_A} \ar[r]^{\tau} &  A \otimes A
                                             \ar[ld]^{m_A} \\
     A & }
 \label{abelian}
\end{equation}

\smallskip

\noindent An important observation is that nonassociative $\mathbb{K}$-algebras
play a key role in the context of Hopf algebras, in particular, for those of combinatorial type. Recall that the Lie algebra $\call(A)$ associated
with a $\mathbb{K}$-algebra $A$ follows from anti-symmetrization of the algebra
product $m_A$. Algebras that give Lie algebras in this way are called Lie
admissible. Another class of Lie admissible algebras are pre-Lie
$\mathbb{K}$-algebras. A {\it left pre-Lie algebra} \cite{Burde_06, Cartier_10,
Manchon_11} is a vector space $V$ equipped with a bilinear product
$\vartriangleright : V \otimes V \to V$ such that the {\it(left) pre-Lie identity},
\begin{equation}
\label{pLrel}
	a \vartriangleright (b\vartriangleright c) - (a\vartriangleright b)
\vartriangleright c
		= b \vartriangleright (a\vartriangleright c) -
(b\vartriangleright a) \vartriangleright c,
\end{equation}
holds for $a,b,c \in V$. This identity rewrites as $L_{[a,b]}=[L_a,L_b]$, where $L_a : V \to V$ is defined by $L_a b:=a \vartriangleright b$. The bracket on the left-hand side is defined by $[a,b]:= a \vartriangleright b - b \vartriangleright a$ and satisfies the Jacobi identity. Right pre-Lie algebras are defined analogously. Note that the (left) pre-Lie identity \eqref{pLrel} can be understood as a relation between associators, i.e., let $\alpha_\vartriangleright : V \otimes V \otimes V \to V$ be defined by
$$
	\alpha_\vartriangleright (a,b,c):=a \vartriangleright (b\vartriangleright c) - (a\vartriangleright b)
\vartriangleright c,
$$
then \eqref{pLrel} simply says that $\alpha_\vartriangleright (a,b,c) = \alpha_\vartriangleright (b,a,c)$. From this it is easy to see that any associative algebra is pre-Lie.

\begin{example} \cite{Burde_06}
Let $A$ be a commutative $\mathbb{K}$-algebra endowed with commuting derivatives $D:=\{ \partial_1,\ldots, \partial_n\}$. For $a \in A$ define $a \partial_i : A \to A$ by  $(a \partial_i)(b):=a\partial_i b$ and $V(n):=\big\{\sum_{i=1}^n a_i \partial_i \ |\ a_i \in A, \partial_i \in D \big\}$. The algebra $(V(n),\lhd)$, where $\sum_{i=1}^n a_i \partial_i \lhd \sum_{j=1}^n a_j \partial_j := \sum_{j,i=1}^n a_j (\partial_ja_i) \partial_i$, is a right pre-Lie algebra called the pre-Lie Witt algebra.
\end{example}

\begin{example} \cite{Burde_06,Cartier_10,Manchon_11}
The next example of a pre-Lie algebra is of a geometric nature and similar to the one above. Let $M$ be a differentiable manifold endowed with a flat and torsion-free connection. The corresponding covariant derivation operator $\nabla$ on the space $\chi(M)$ of vector fields on $M$ provides it with a left pre-Lie algebra structure, which is defined via $a \vartriangleright b:=\nabla_ab$ by virtue of the two equalities $\nabla_ab-\nabla_ba=[a,b]$, $\nabla_{[a,b]}=[\nabla_a,\nabla_b]$. They express the vanishing of torsion and curvature respectively. Let $M=\R^n$ with its standard flat connection. For $a=\sum_{i=1}^na_i\partial_i$ and $b=\sum_{i=1}^nb_i\partial_i$ it follows that
\begin{equation*}
	a\vartriangleright b=\sum_{i=1}^n\left(\sum_{j=1}^na_j(\partial_jb_i)\right)\partial_i.
\end{equation*}
\end{example}

\begin{example} \cite{Burde_06}
Denote by $W$ the space of all words over the alphabet $\{a,b\}$. Define for words $v$ and $w=w_1 \cdots w_n$ in $W$ the product $w \circ v := \sum_{i=0}^{n} \epsilon(i) w \lhd_i v$, where $w \lhd_i v$ denotes inserting the word $v$ between letters $w_i$ and $w_{i+1}$ of $w$, i.e., $w \lhd_i v = w_1 \cdots w_i v w_{i+1} \cdots w_n$ and
$$
\epsilon(i):=
\begin{cases}
 -1, & w_i=a, w_{i+1}=b\\
+1, & w_i=b, w_{i+1}=a\ or\ \emptyset\\
+1, & w_i=\emptyset, w_{i+1}=a\\
 0,  & \mbox{otherwise}.
\end{cases}
$$
The algebra $(W,\circ)$ is right pre-Lie. For example,
$$
	a \circ a = aa,\quad
	a \circ ab = aba,\quad
	ab \circ a = aab - aab + aba = aba, \quad
	ba \circ ab = baba.
$$
\end{example}


\subsubsection{Coalgebra}
\label{ssubsec:coalgebras}

The definition of a $\mathbb{K}$-coalgebra is most easily obtained by
reversing the arrows in diagrams (\ref{assoc2}) and (\ref{unit}). Thus, a
$\mathbb{K}$-{\it{coalgebra}} is a triple $(C,\Delta_{C},\epsilon_{C})$, where
$C$ is  a $\mathbb{K}$-vector space carrying a {\it{coassociative coproduct}} map
$\Delta_C: C \to C \otimes C$, i.e., $(\Delta_C \otimes \id_C) \circ
\Delta_C=(\id_C \otimes \Delta_C) \circ \Delta_C : C \to C \otimes C \otimes C$,
and $\epsilon_C: C \to \mathbb{K}$ is the counit map which satisfies
$(\epsilon_C \otimes \id_C) \circ \Delta_C = \id_C  =(\id_C \otimes \epsilon_C)
\circ \Delta_C$. Its kernel $\mathrm{ker}(\epsilon_C) \subset C$ is called the
{\it{augmentation ideal}}.

A simple example of a coalgebra is the field $\mathbb{K}$ itself with the
coproduct $ \Delta_\mathbb{K}: \mathbb{K} \to \mathbb{K} \otimes \mathbb{K}$, $c
\mapsto c \otimes 1$ and $\epsilon_\mathbb{K} := \id_\mathbb{K}: \mathbb{K} \to
\mathbb{K}$.

Using Sweedler's notation for the coproduct of an element $x \in C$, $\Delta_C(x) = \sum_{(x)} x^{(1)}
\otimes x^{(2)}$, provides a simple description of coassociativity
$$
  	\sum_{(x)}\bigg(\sum_{(x^{(1)})} x^{(1)(1)} \otimes x^{(1)(2)}\bigg ) \otimes x_{(2)}
 	= \sum_{(x)} x^{(1)} \otimes \bigg (\sum_{(x^{(2)})} x^{(2)(1)}\otimes x^{(2)(2)} \bigg).
$$
It permits the use of a transparent notation for iterated coproducts: $\Delta_C^{(n)}:C \to C^{\otimes n+1}$, where $\Delta^{(0)}:=\id$,  $\Delta_C^{(1)}:=\Delta_C$, and
$$
	\Delta_C^{(n)}:=(\id \otimes \Delta_C^{(n-1)}) \circ \Delta_C=(\Delta_C^{(n-1)} \otimes \id_{C}) \circ \Delta_C.
$$
For example,
$$
    	\Delta_C^{(2)}(x):=(\id_C \otimes \Delta_C) \circ \Delta_C(x)
        = (\Delta_C \otimes \id_C) \circ \Delta_C(x) = \sum_{(x)} x^{(1)} \otimes x^{(2)} \otimes x^{(3)}.
$$
A {\it cocommutative} coalgebra satisfies $\tau \circ \Delta_{C}=\Delta_{C}$, which amounts to reversing the arrows in diagram (\ref{abelian}). An element $x$ in a coalgebra $(C,\Delta_C,\epsilon_C)$ is called {\it{primitive}} if
$\Delta_C(x) = x \otimes 1_C + 1_C \otimes x$. It is called {\it group-like} if $\Delta_C(x) = x \otimes x$. The set of primitive elements in $C$ is denoted by $P(C)$. A subspace $I \subseteq C$ of a $\mathbb{K}$-coalgebra $(C,\Delta_C,\epsilon_C)$ is a {\it{subcoalgebra}} if $\Delta_C(I) \subseteq I \otimes I$. A subspace $I \subseteq C$ is called a (left-, right-) {\it{coideal}} if ($\Delta_C(I) \subseteq I \otimes C$, $\Delta_C(I) \subseteq  C \otimes I$) $\Delta_C(I) \subseteq I \otimes C + C \otimes I$.


\subsubsection{Bialgebra}
\label{ssubsec:bialgebras}

A $\mathbb{K}$-{\it{bialgebra}} consists of a $\mathbb{K}$-algebra and a $\mathbb{K}$-coalgebra which are compatible \cite{Abe_80,Figueroa-Gracia-Bondia_05,Manchon_01,Radford_12,Sweedler_69}. More precisely, a $\mathbb{K}$-bialgebra is a quintuple $(B,m_B,\eta_B,\Delta_B,\epsilon_B)$, where $(B,m_B,\eta_B)$ is a $\mathbb{K}$-algebra, and $(B,\Delta_B,\epsilon_B)$ is a $\mathbb{K}$-coalgebra, such that $m_B$ and $\eta_B$ are morphisms of $\mathbb{K}$-coalgebras with the natural coalgebra structure on the space $B \otimes B$. Commutativity of the following diagrams encodes the compatibilities
\begin{equation}
\xymatrix{
          B \otimes B \ar[rr]^{m_B} \ar[d]_{\scriptstyle{\tau_2
          (\Delta_B \otimes \Delta_B)}}  & & B \ar[d]^{\Delta_B}\\
          B \otimes B \otimes B \otimes B \ar[rr]^{\qquad m_B \otimes m_B} && B
\otimes B }
    \label{compat1}
\qquad
  \xymatrix{
            B \otimes B \ar[rr]^{\epsilon_B \:\otimes \epsilon_B} \ar[d]_{m_B} &
&
            \mathbb{K} \otimes \mathbb{K} \ar[d]^{m_\mathbb{K}}\\
            B \ar[rr]^{\epsilon_B} && \mathbb{K} }
\end{equation}

\begin{equation}
 \xymatrix{
             \mathbb{K} \ar[rr]^{\eta_B} \ar[d]_{\Delta_\mathbb{K}} && B
\ar[d]_{\Delta_B}\\
             \mathbb{K} \otimes \mathbb{K} \ar[rr]^{\eta_B \:\otimes \eta_B} &&
B \otimes B
            }
\qquad
  \xymatrix{
          \mathbb{K} \ar[rr]^{\eta_B} \ar[rd]_{\id_\mathbb{K}} && B
\ar[ld]^{\epsilon_B} \\
           & \mathbb{K} & },
\end{equation}
where $\tau_2:=(\id_B \otimes \tau \otimes \id_B)$. Equivalently, $\Delta_B$ and
$\epsilon_B$ are morphisms of $\mathbb{K}$-algebras, with the natural algebra
structure on the space  $B\otimes B$. By a slight abuse of notation one writes
$\Delta_B(m_B(b \otimes b')) = \Delta_B(b)\Delta_B(b')$ for $b,b' \in B$,
saying that the {\it{coproduct of the product is the product of the coproduct}}.
The identity element in $B$ will be denoted by $\mathbf{1}_B$, and all algebra
morphisms are required to be unital. Note that if $x_1,x_2$ are primitive in
$B$, then $[x_1,x_2]:=m_B(x_1 \otimes x_2) - m_B( x_2 \otimes x_1)$ is primitive
as well, i.e., the set $P(B)$ of primitive elements of a bialgebra $B$ is a Lie
subalgebra of the Lie algebra $\call(B)$.

\smallskip

A bialgebra $B$ is called {\it{graded}} if there are $\mathbb{K}$-vector
spaces $B_{n}$, $n \geq 0$, such that
\begin{enumerate}
	\item $B= \bigoplus_{n \geq 0} B_{n}\ ,$
	
	\item $m_B(B_{n} \otimes B_{m}) \subseteq B_{n+m},$
	
	\item $\Delta_B(B_{n}) \subseteq \bigoplus_{p+q=n} B_{p} \otimes B_{q}.$
\end{enumerate}
Elements $x \in B_{n}$ are given a degree $\mathrm{deg}(x)=n$. For a
{\it{connected}} graded bialgebra $B$, the degree zero part is $B_{0} =
\mathbb{K}\mathbf{1}_B$. Note that $\mathbf{1}_B$ is group-like. A graded
bialgebra $B = \bigoplus_{n \geq 0} B_{n}$ is said to be of {\it{finite type}}
if its homogeneous components $B_{n}$ are $\mathbb{K}$-vector spaces of finite
dimension.

Let $B$ be a connected graded $\mathbb{K}$-bialgebra. One can show \cite{Manchon_01} that the
coproduct of any element $x \in B$ is given by
$$
    \Delta_B(x) = x \otimes \mathbf{1}_B + \mathbf{1}_B \otimes x + \sideset{}{'}\sum_{(x)} x' \otimes x'',
$$
where
$$
	\Delta'(x):=\sideset{}{'}\sum_{(x)}  x{'} \otimes x{''} \in
	\bigoplus_{p+q=n \atop p>0,\ q>0} B_{p} \otimes B_{q}
$$
is the {\it reduced coproduct}, which is coassociative on the augmentation ideal
$\mathrm{ker}(\epsilon_B) := B^+:= \bigoplus_{n>0} B_{n}$. Elements in the kernel of
$\Delta'_B$ are primitive elements of $B$.

\begin{example} \label{exam:dividedpowers}
{\it{Divided powers}} $(D,m_D,\eta_D,\Delta_D,\epsilon_D)$ are a graded
bialgebra, where $D=\bigoplus_{n=0}^\infty D_n$, $D_n:= \mathbb{K} d_n$. The
product is given by $m_D (d_m \otimes d_n)= {m+n \choose m} d_{m+n}$, and the
unital map is $\eta_D: \mathbb{K} \to D$, $1_{\mathbb{K}} \mapsto
d_{0}:=\mathbf{1}_D$. The coproduct $\Delta_D: D \to D \otimes D$ maps $d_{n}
\mapsto \sum_{k=0}^{n} d_{k} \otimes d_{n-k}$, and $\epsilon_D: D \to
\mathbb{K}$, $d_{n} \mapsto \delta_{0,n}1_{\mathbb{K}}$, where $\delta_{0,n}$ is
the usual Kronecker delta.
\end{example}

\medskip

For a $\mathbb{K}$-algebra $A$ and a $\mathbb{K}$-coalgebra $C$, the
{\it{convolution product}} of two linear maps $f,g  \in
L(C,A):={\rm{Hom}}_\mathbb{K}(C,A)$ is defined to be the linear map $f \star
g \in L(C,A)$ given for $a \in C$ by
\begin{equation}
  \label{def:convolution}
    (f \star g)(a) := m_{A} \circ (f \otimes g)\circ \Delta_C(a) = \sum_{(a)} f(a^{(1)}) \, g(a^{(2)}).
\end{equation}
In other words
$$
  C \xrightarrow{\Delta_C}  C \otimes C \xrightarrow{f \otimes g} A \otimes A
\xrightarrow{m_A} A.
$$
It is easy to see that associativity of $A$ and coassociativity of $C$ imply
the following.

\begin{theorem} {\cite{Abe_80,Figueroa-Gracia-Bondia_05,Manchon_01,Radford_12,Sweedler_69}}
$L(C,A)$ with the convolution product \eqref{def:convolution} is an unital
associative $\mathbb{K}$-algebra with unit $\eta:=\eta_A \circ \epsilon_C$.
\end{theorem}

The algebra $A$ can be replaced by the base field $\mathbb{K}$. For a bialgebra
$B$ the theorem describes the convolution algebra structure on $L(B,B)$ with
unit $\eta:=\eta_B \circ \epsilon_{B}$.

For the maps $f_i \in L(C,A)$, $i=1,\ldots, n$, $n>1$, multiple convolution
products are defined by
\begin{equation}
  \label{def:convpowers}
    f_1 \star f_2 \star \dots \star f_n :=
        m_{A} \circ (f_1 \otimes f_2 \otimes \dots \otimes f_n) \circ
\Delta_C^{(n-1)}.
\end{equation}
Recall that $\Delta_C^{(0)} := \id_C$, and for $n>0$, $\Delta_C^{(n)}:=(\Delta_C^{(n-1)} \otimes \id_{C}) \circ \Delta_C$.


\subsubsection{Hopf algebra}
\label{ssubsec:hopfalgebras}

\begin{definition}\label{def:HopfAlgebra}  {\cite{Abe_80,Radford_12,Sweedler_69}}
A {\it{Hopf algebra}} is a $\mathbb{K}$-bialgebra
$(H,m_H,\eta_H,\Delta_H,\epsilon_H,S)$ together with a particular
$\mathbb{K}$-linear map $S: H \to H$ called the {\it{antipode}}, which satisfies the
Hopf algebra axioms \cite{Abe_80,Radford_12,Sweedler_69}. The algebra unit in $H$ is
denoted by $\mathbf{1}_H$.
\end{definition}

The antipode map has the property of being an antihomomorphism for both the algebra
and the coalgebra structures, i.e., $S(m_H(x \otimes y))=m_H(S(y) \otimes S(x))$
and $\Delta_H \circ S = (S \otimes S)\circ \tau \circ \Delta_H$. The
necessarily unique antipode $S \in L(H,H)$ is the inverse of the identity map
$\id_H : H \to H$ with respect to the convolution product
\begin{equation}
     S \star \id_H = m_H  \circ (S \otimes \id_H)\circ \Delta_H
                  = \eta_H \circ \epsilon_H
                  = m_H \circ (\id_H \otimes S)\circ \Delta_H= \id_H \star S .
    \label{eq:antipode}
\end{equation}
If the Hopf algebra $H$ is commutative or cocommutative, then $S \circ S=\id_H$.

Recall that the universal enveloping algebra $\mathcal{U}(\call)$ of a Lie algebra $\call$ has the structure of a Hopf algebra \cite{Reutenauer_93} and provides a natural example. An important observation is contained in the next result.

\begin{proposition}\label{prop:conngradedbialg}  {\cite{Figueroa-Gracia-Bondia_05}}
Any connected graded bialgebra $H=\bigoplus_{n \ge 0} H_n$ is a connected graded Hopf algebra. The antipode is defined by the geometric series $S:=\id_H^{\star (-1)} = \big(\eta_H \circ \epsilon_H - (\eta_H
\circ \epsilon_H - \id_H)\big)^{\star (-1)}$.
\end{proposition}

See \cite{Figueroa-Gracia-Bondia_05} for a proof and more details. Hence, for any $x \in H_n$ the antipode is computed by
\begin{equation}
\label{def:antipode1}
	S(x) = \sum_{k\ge 0}\big(\eta_H \circ \epsilon_H - \id_H\big)^{\star k}(x).
\end{equation}
It is well-defined as the sum on the right-hand side terminates at $\mathrm{deg}(x)=n$
due to the fact that the projector $P := \id_H - \eta_H \circ \epsilon_H$ maps $H$ to its augmentation ideal $\mathrm{ker}(\epsilon_H)$. Note that the antipode preserves the grading, i.e., $S(H_{n}) \subseteq H_{n}$.

\begin{corollary} \label{cor:recursiveS}  {\cite{Figueroa-Gracia-Bondia_05}}
The antipode $S$ for a connected graded Hopf algebra $H = \bigoplus_{n\ge0}
H_{n}$ may be defined recursively in terms of either of the two formulae
\begin{subequations} \label{antipode2}
\begin{align}
        S(x) &= - S \star P(x) = -x - \sideset{}{'}\sum_{(x)} S(x{'})x{''}, \label{eq:left-antipode-recursion} \\
        S(x) &= - P \star S (x) = -x - \sideset{}{'}\sum_{(x)} x{'}S(x{''}),
\end{align}
\end{subequations}
for $x \in \mathrm{ker}(\epsilon_H)= \bigoplus_{n>0} H_{n}$, which follow
readily from (\ref{eq:antipode}) and $S(\mathbf{1}_H)=\mathbf{1}_H$.
\end{corollary}

{{These recursions make sense due to the fact that on the righthand side the antipode is calculated on elements $x'$ or $x''$, which are of strictly smaller degree than the element $x$.}} Let $A$ be a $\mathbb{K}$-algebra and $H$ a Hopf algebra. An element $\phi \in L(H,A)$ is called a {\it{character}} if $\phi$ is a unital algebra morphism, that is, $\phi(\mathbf{1}_H) =\mathbf{1}_A$ and
\begin{equation}
\label{character}
	\phi \big( m_H(x \otimes y) \big) = m_A\big(\phi(x) \otimes
\phi(y)\big).
\end{equation}
An {\it{infinitesimal character}} with values in $A$ is a linear map $\xi \in
L(H, A)$ such that for $x, y \in \mathrm{ker}(\epsilon_H)$, $\xi(m_H(x \otimes
y)) = 0$, which implies $\xi(\mathbf{1}_H) =0$. An equivalent way to
characterize  infinitesimal characters is as {\it{derivations}}, i.e.,
\begin{equation}
 \label{def:inifiniChar}
    \xi\big(m_H(x \otimes y)\big)
    = m_A \big(\eta_A \circ \epsilon_H(x)\otimes \xi(y)\big) + m_A \big(\xi(x)\otimes \eta_A \circ
\epsilon_H(y)\big)
\end{equation}
for any $x,y \in H$. The set of characters (respectively  infinitesimal characters)
is denoted by $G_A \subset L(H,A)$ (respectively $g_A \subset L(H,A)$).

Let $A$ be a commutative $\mathbb{K}$-algebra. The linear space of infinitesimal
characters, $g_{A}$, forms a Lie algebra with respect to the Lie bracket defined on $L(H,A)$
in terms of the convolution product
$$
	[\alpha , \beta] := \alpha \star \beta - \beta \star \alpha.
$$
Moreover, $A$-valued characters, $G_A$, form a group. The inverse is given by
composition with the antipode $S$ of $H$. Since $\alpha(\mathbf{1}_H) =0$ for $\alpha \in g_{A}$, the exponential defined by its
power series with respect to convolution, $\exp^\star(\alpha)(x) := \sum_{j\ge
0} \frac{1}{j!}\alpha^{\star j}(x)$, is a finite sum terminating at $j=n$ for any $x \in H_n$.

\begin{proposition}\label{prop:exp} {\cite{Figueroa-Gracia-Bondia_05,Manchon_01}}
$\exp^\star$ restricts to a bijection from $g_A$ onto $G_A$.
\end{proposition}

The compositional inverse of $\exp^\star$ is given by the logarithm defined with respect to the convolution product, $\log^\star(\eta_A \circ \epsilon_H +\gamma)(x)=\sum_{k\ge 1}{(-1)^{k-1}\over k}\gamma^{\star k}(x)$, where $\gamma \in g_A$. Again the sum terminates at $k=n$ for any $x \in H_n$ as $\gamma(\mathbf{1}_H) =0$. For details and proofs of these well-known facts the reader is referred to \cite{Figueroa-Gracia-Bondia_05,Manchon_01}.

\begin{remark} \label{rmk:milnormoore} An important result which is due to {{Milnor and Moore}} \cite{Milnor-Moore_65} concerns the structure of cocommutative connected graded Hopf algebras of finite type. It states that any such Hopf algebra $H$ is isomorphic to the universal enveloping algebra of its primitive elements, i.e., $H \cong \mathcal{U}(P(H))$. See also \cite{FGV01}.
\end{remark}

\begin{remark} \label{rmk:groupscheme}
An observation concerning the relationship between the group $G_A
\subset L(H, A)$ and the Hopf algebra $H$ will be important in
the analysis which follows. The reader is referred to the paper of Manchon and
Frabetti \cite{Frabetti-Manchon_15} for details and additional references. By
definition, elements in $G_A$ map all of $H$ into the commutative unital algebra
$A$. However, an element $x \in H$ can also be seen as an $A$-valued function on
$G_A$. Indeed, let $\Phi \in G_A$, then $x(\Phi):=\Phi(x) \in A$ and the usual
pointwise product of functions
$
	(xy)(\Phi)=x(\Phi)y(\Phi)
$
follows from \eqref{character} since $\Phi \in G_A$. The definition of the
convolution product \eqref{def:convolution} in terms of the coproduct of $H$
implies a natural coproduct on functions $x \in H$, that is,
$\Delta(x)(\Phi,\Psi) :=(\Phi \star \Psi)(x) \in A$. Similarly, the inverse of
$G_A$ as well as its unit correspond naturally to the antipode and counit map on
$H$, respectively. This {\em reversed} perspective on the relationship between
$H$ and its group of characters $G_A$ allows one to interpret $H$ as the (Hopf)
algebra of  {\it{coordinate functions}} of the group $G_A$. More precisely, $H$
contains the {\em representative functions} over $G_A$. The reader is directed
to Cartier's work \cite{Cartier_07} for a comprehensive review of this topic. In
the context of input-output systems in nonlinear control theory, a particular
group of unital Fliess operators is central. Its product, unit and inverse are used to identify its Hopf algebra of coordinate functions.
\end{remark}

\medskip

\begin{example} \label{exam:shuffle}
Three examples of Hopf algebra are presented: the unshuffle, shuffle {\cite{Berstel-Reutenauer_88}} and
quasi-shuffle Hopf algebras {\cite{Hoffman_00}}.

\begin{enumerate}

\item \label{concat}
{{
Let $X:=\{{x_1},{x_2}, {x_3},\ldots,x_m\}$ be an alphabet with $m$ letters. As before, $X^*$ is the set of words with letters in $X$. The length of a word $\eta=x_{i_1}\cdots x_{i_n}$ in $X^*$ is defined by the number of letters it contains, i.e., $|\eta|:=n$. The
empty word $\mathbf{1} \in X^*$ has length zero. The vector space ${\mathbb K}\langle X\rangle$, which is freely generated by
$X^*$ and graded by length, becomes a noncommutative, unital, connected, graded
algebra by concatenating words, i.e., for $\eta=x_{i_1}\cdots x_{i_n}$ and
$\eta'=x_{j_1}\cdots x_{j_l}$, $\eta \cdot \eta':= x_{i_1}\cdots x_{i_n}
x_{j_1}\cdots x_{j_l}$, and $|\eta \cdot \eta'| = n+l$. The {\it unshuffle
coproduct} is defined by declaring the elements in $X$ to be primitive, i.e.,
$\Delta^{\!\!\shuffle}(x_i):= {x_i} \otimes \mathbf{1} + \mathbf{1} \otimes
{x_i}$ for all $x_i \in X$, and by extending it multiplicatively. In this case, ${\mathbb
K}\langle X\rangle$ is turned into the cocommutative {\it unshuffle Hopf
algebra} $H_{conc}$. For instance, for the letters $x_{i_1},x_{i_2} \in X$ the coproduct of the length two word $\eta=x_{i_1}x_{i_2}$ is
$$
	\Delta^{\!\!\shuffle}(\eta) =
\Delta^{\!\!\shuffle}(x_{i_1})\Delta^{\!\!\shuffle}(x_{i_2})
	= x_{i_1}x_{i_2} \otimes \mathbf{1}  + \mathbf{1}  \otimes x_{i_1}x_{i_2}
	+ x_{i_1} \otimes x_{i_2} + x_{i_2} \otimes x_{i_1}.
$$
The general form of $\Delta^{\!\!\shuffle}$ for an arbitrary word $\eta=x_{i_1}\cdots x_{i_l} \in X^*$ is given by
\begin{align}
\label{unshuffle}
	\Delta^{\!\!\shuffle}(\eta)
	=\prod_{j=1}^{l}\Delta^{\!\!\shuffle}(x_{i_j})=\sum_{\alpha, \beta \in X^*}
	\langle \alpha \shuffle  \beta, {\eta} \rangle \alpha \otimes  \beta.
\end{align}
The coefficient in the sum over words $\alpha, \beta \in X^*$ on the right-hand
side is defined through the linearly extended bracket $\langle \alpha , \nu
\rangle : = 1$ if ${\alpha}={\nu}$, and zero otherwise. The product $\shuffle$
displayed in \eqref{unshuffle} is the shuffle product of words
\eqref{shuffleproduct} introduced above. The antipode of $H_{conc}$ turns out to be
\begin{equation}
\label{antipode1}
	S(x_{i_1}\cdots x_{i_l})=(-1)^l x_{i_l}\cdots x_{i_1}.
\end{equation}
It is interesting to check that \eqref{antipode1} satisfies the recursions \eqref{antipode2}.
}}
{{
\item \label{shuf}
The same space ${\mathbb K}\langle X\rangle$ can be turned into a unital,
connected, graded, commutative, noncocommutative Hopf algebra, $H_{\shuffle}$, known
as {\it shuffle Hopf algebra}, by defining its algebra structure in terms of the
shuffle product on words \eqref{shuffleproduct}, and its coproduct by {\it
deconcatenation}. The latter is defined on words $\eta=x_{i_1}\cdots x_{i_l} \in
X^*$ as follows
\begin{align}
\label{decon}
	\Delta(\eta)=\eta \otimes \mathbf{1}  + \mathbf{1}  \otimes \eta +
	\sum_{k=1}^{l-1} x_{i_1}\cdots x_{i_k} \otimes x_{i_{k+1}}\cdots x_{i_l}.
\end{align}
It is easy to show -- and left to the reader, -- that this gives a coassociative
coproduct which is compatible with the shuffle product. The antipode of
$H_{\shuffle}$ is the same as that of $H_{conc}$, i.e., $S(x_{i_1}\cdots x_{i_l}):=(-1)^l x_{i_l}\cdots x_{i_1}$.
Again, it is interesting to  verify that  it satisfies both recursions \eqref{antipode2}.
}}
{{
\item \label{quasishuf}
The last example of a connected graded Hopf algebra is defined on the countable alphabet $A$. Here $A^*$ denotes the monoid of words $w=a_{i_1} \cdots a_{i_l}$ generated by the letters from $A$ with concatenation as product. The degree of an element $a_i
\in A$ is defined to be $\deg(a_i)$. It is extended to words additively, i.e., $\deg(a_{i_1} \cdots a_{i_l})= \deg(a_{i_1}) + \cdots + \deg(a_{i_l})$. The empty word $\mathbf{1} \in A^*$ is of degree zero. Moreover, it is assumed that $A$ itself is a graded commutative semigroup with bilinear product $[- \ -] \colon  A \times A \to A$. The degree $\deg([a_{i}\ a_{j}]):=\deg(a_{i}) +\deg(a_{j})$. Commutativity and associativity of $[- \ -]$ allow for a notational simplification, i.e., $[a_{i_1} \cdots  a_{i_n}] :=[a_{i_1}\ [\cdots [a_{i_{n-1}}\ a_{i_{n}}]]\cdots]$. The free noncommutative algebra of words $w=a_{i_1} \cdots a_{i_l}$ over the alphabet $A$ is denoted ${\mathbb K}\langle A\rangle$. The commutative and associative quasi-shuffle product on words $w,v \in {\mathbb K}\langle A\rangle$ is defined by
\begin{enumerate}
 	\item[i)]\ $\mathbf{1} \star v :=v \star \mathbf{1} := v$,
 	\item[ii)]\ $a_i v \star a_j w := a_i(v \star a_j w) +a_ j(a_i v \star w) + [a_i \ a_j] (v \star w)$,
\end{enumerate}
where $a_i,a_j$ are letters in $A$. For instance,
\begin{eqnarray*}
	a_{i_1} \star a_{i_2} &=& a_{i_1}a_{i_2} + a_{i_2}a_{i_1} + [a_{i_1}\ a_{i_2}]\\
	a_{i_1} \star  a_{i_2}a_{i_3} &=& a_{i_1} a_{i_2}a_{i_3} + a_{i_2}a_{i_1}a_{i_3} + a_{i_2}a_{i_3}a_{i_1}
					+ [a_{i_1}\ a_{i_2}]a_{i_3} + a_{i_2}[a_{i_1}a_{i_3}].
\end{eqnarray*}
Hoffman \cite{Hoffman_00} showed that the quasi-shuffle algebra $H_\star$ is a Hopf algebra with respect to the deconcatenation coproduct \eqref{decon}. The antipode $S:H_{\star} \to H_{\star}$ is deduced from the recursion \eqref{eq:left-antipode-recursion}, e.g.,
\begin{eqnarray}
        \label{antipodeQuasi}
        S(a_{i_1} \cdots  a_{i_n}) &=& - (S \star P)(a_{i_1} \cdots  a_{i_n}) \nonumber \\
        				   &=& \big(- \id_{H_{\star}} - m_{\star}
\circ(S \otimes  \id_{H_{\star}}) \circ \Delta' \big)(a_{i_1} \cdots  a_{i_n})\\
        				   &=& -a_{i_1} \cdots  a_{i_n} -
\sum_{l=1}^{n-1} S(a_{i_1} \cdots  a_{i_l}) \star a_{i_{l+1}} \cdots  a_{i_n},
\end{eqnarray}
where the projector $P := \id_{H_{\star}} - \eta_{H_{\star}} \circ
\epsilon_{H_{\star}}$ maps $H_{\star}$ to its augmentation ideal $\mathrm{ker}(\epsilon_{H_{\star}})$. For example, the antipode for the letter $a_i$ and the word $a_{i}a_{j}$ are respectively
$$
	S(a_i)=-a_i,
	\quad
	S(a_{i}a_{j})=-a_{i}a_{j} + a_{i}\star a_{j}  = a_{j}a_{i} + [a_{i}\ a_{j}].
$$
If $[a_i\ a_j]=0$ for any letters $a_i,a_j \in A$, then the quasi-shuffle product reduces to the ordinary shuffle product \eqref{shuffleproduct} on words, and $H_\star$ reduces to  $H_\shuffle$.
}}

\end{enumerate}
\end{example}
In Section \ref{sec:cha-crt} another example of a connected graded Hopf algebra is presented, one which plays a key role in the context of the output feedback interconnection. In Section \ref{sec:towards-discretization} the quasi-shuffle product is employed in the context of products of discrete-time Fliess operators.


\section{Algebras Induced by the Interconnection of Fliess Operators}
\label{sec:system-interconnections}

In engineering applications, where input-output systems are represented in terms of Fliess operators, it is natural to interconnect systems to create models of more complex systems. It is known that all the basic interconnection types, such as the parallel, cascade and feedback connections, are well-posed. For example, under suitable assumptions, the output of a given Fliess operator generates an admissible input for another Fliess operator in the cascade connection. In each case it is also known that the aggregate system has a Fliess operator representation. Therefore, a family of formal power series products is naturally induced whereby the generating series of an aggregate system can be computed in terms of the generating series of its subsystems. Of particular interest here are the cascade and feedback interconnections, as their respective
products define a semi-group and group that are central in control theory. In this section, these algebraic structures are described in detail.


\subsection{Cascade interconnections}
\label{subsec:nonrecursive-interconnections}

\begin{figure}[t]
\begin{center}
\includegraphics[scale=0.38]{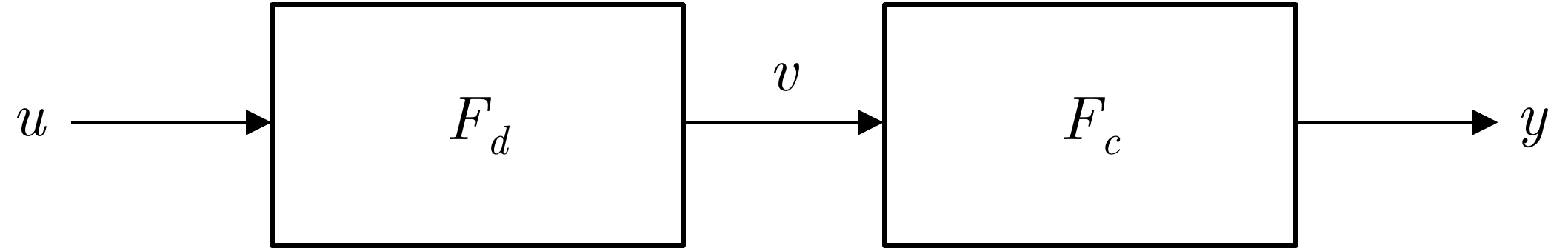} \\[0.15in]
\includegraphics[scale=0.38]{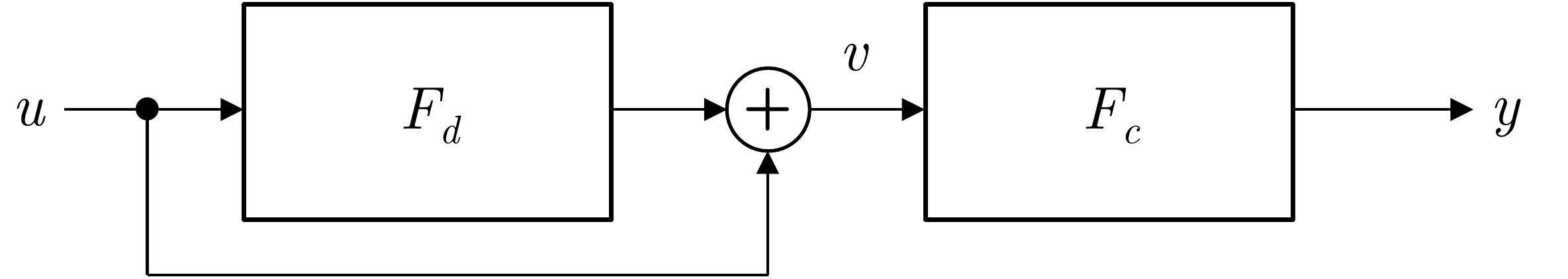}
\end{center}
\caption{Fliess operator cascades yielding the composition product (top) and
modified composition product (bottom)}
\label{fig:composition-connections}
\end{figure}

Consider the cascade interconnections of two Fliess operators shown in
Figure~\ref{fig:composition-connections}.
The first is a simple cascade interconnection corresponding to a direct
composition of the operators $F_c$ and $F_d$,
namely, $F_c \circ F_d$. The other involves a {\em direct feed term} passing
from the input $u$ to the input $v$ so that $F_c\circ (I+F_d)$, where $I$
denotes the identity operator.  This direct feed term is a common feature found
in some
control systems and is particularly important in feedback systems as will be
discussed shortly. The primary claim is that each cascade interconnection
induces a locally finite product on the level of formal power series,
which unambiguously describes the interconnected system as
generating series of Fliess operators are known to be unique
\cite{Fliess_83,Wang_90}. Specifically, the {\em composition product} satisfies
$F_c \circ F_d=F_{c\circ d}$, and the {\em modified composition product}
satisfies $F_c \circ (I+F_d) = F_{c\modcomp d}$.
Each product is defined in terms of a certain algebra homomorphism. The morphism
for the modified composition product ultimately defines a pre-Lie product, which
is at the root of all the underlying combinatorial structures at play.

For a fixed $d \in \allseriesm$ and alphabet $X=\{x_0,x_1,\ldots,x_m\}$, let
$\psi_d$ be the continuous (in the ultrametric sense) algebra homomorphism
mapping $\allseries$ to the set of vector space endomorphisms $\Endallseries$
uniquely specified by
$\psi_d(x_i\eta)=\psi_d(x_i)\circ \psi_d(\eta)$ for $x_i\in X$, $\eta\in X^\ast$ with
\begdi \label{eq:psi-d-on-words}
	\psi_d(x_i)(e)=x_0(d_i\shuffle e),
\enddi
$i=0,1,\ldots,m$ and any $e\in\allseries$, and where $d_i$ is the $i$-th
component series of $d \in \allseriesm$ ($d_0:=\emptyset$). By definition,
$\psi_d(\emptyset)$ is the identity map on $\allseries$. The following theorem
describes the generating series for the direct cascade connection of two Fliess
operators.

\begth {\rm {\cite{Ferfera_79,Ferfera_80,Gray-Li_05,Thitsa-Gray_12}}}
Given any {$c\in\allseriesellLC$ and $d\in\allseriesmLC$}, the composition $F_c\circ
F_d=F_{c\circ d}$, where the \bfem{composition product} of $c$ and $d$ is given
by
\begdi
	c\circ d=\sum_{\eta\in X^\ast} (c,\eta)\,\psi_d(\eta)(\mathbf{1}),
\enddi
{and $c\circ d\in\allseriesellLC$.}
\endth

It is not difficult to show that this composition product is locally finite (and
therefore summable), associative and $\re$-linear in its left argument. As this
product lacks an identity, it defines a semi-group on $\allseriesm$. In
addition, this product distributes to the left over the shuffle product, which
reflects the fact that in general
$$
	F_{(c\shuffle d)\circ e}=(F_cF_d) \circ F_e = F_c[F_e]F_d[F_e]
	=F_{(c\circ e)\shuffle(d\circ e)}.
$$
Finally, it is known that the composition product defines an ultrametric
contraction on $\allseriesm$.

\begin{example} \label{ex:LTI-cascade-detail}
In the case of two linear time-invariant systems with analytic kernels
$h_c(t)=\sum_{i\geq 0} (c,x_0^ix_1)t^i/i!$
and $h_d(t)=\sum_{i\geq 0} (d,x_0^ix_1)t^i/i!$, respectively, a direct
calculation gives the kernel for the composition
\begin{align*}
	(h_c\ast h_d)(t)	&:=\int_0^t h_c(t-\tau)h_d(\tau)\,d\tau
				=\sum_{k=1}^\infty \left[\sum_{j=0}^{k-1} (c,x_0^{k-j-1}x_1)(d,x_0^jx_1)\right] \frac{t^k}{k!} \\
				&= \sum_{k=1}^\infty (c\circ d,x_0^kx_1)
\frac{t^k}{k!}=:h_{c\circ d}(t).
\end{align*}
\end{example}

The introduction of a direct feed term in the composition interconnection
requires a modification to the set up, namely, the new algebra homomorphism
$\phi_d$ from $\allseries$ to $\Endallseries$ where
$\phi_d(x_i\eta)=\phi_d(x_i)\circ \phi_d(\eta)$ for $x_i\in X$, $\eta\in X^\ast$ with
\begdi \label{eq:phi-d-on-words}
	\phi_d(x_i)(e)=x_ie+x_0(d_i\shuffle e),
\enddi
$i=0,1,\ldots,m$ and any $e\in\allseries$, and where $d_0:=0$. Again,
$\phi_d(\emptyset)$ is the identity map on $\allseries$. The direct feed
term is encoded in the new term $x_ie$ shown above. This yields the desired
formal power series product as described next.

\begth {\rm \cite{Gray-Li_05,Li_04}}
Given any {$c\in\allseriesellLC$ and $d\in\allseriesmLC$,} the composition $F_c\circ
(I+F_d)=F_{c\modcomp d}$, where
the \bfem{modified composition product} of $c$ and $d$ is given by
\begdi
	c\modcomp d=\sum_{\eta\in X^\ast} (c,\eta)\, \phi_d(\eta)(\mathbf{1}),
\enddi
{and $c\modcomp d\in\allseriesellLC$.}
\endth

This product is also always summable, but it is not associative, in fact, in
general
\begeq
	(c\modcomp d)\modcomp e=c\modcomp(d\modcomp e+e)
\label{eq:cmod-non-associative}
\endeq
\cite{Li_04}. The following lemma describes some other elementary properties of
the modified composition product.

\begle \label{le:modified-product-properties} {\rm
\cite{Gray-etal_SCL14,Gray-Li_05}}
The modified composition product
\begin{description}
\item[(1)] is left $\re$-linear;
\item[(2)] satisfies $c\modcomp 0=c$;
\item[(3)] satisfies $c\modcomp d=k\in\re^\ell$ for any fixed $d$ if and only if
$c=k$;
\item[(4)] satisfies $(x_0c)\modcomp d=x_0(c\modcomp d)$ and
	$(x_ic)\modcomp d=x_i(c\modcomp d)+x_0(d_i\shuffle (c\modcomp d))$;
\item[(5)] distributes to the left over the shuffle product.
\end{description}
\endle

\begin{table}[t]
\caption{Composition products involving $c$, $d$, $c_\delta=\delta+c$,
$d_\delta=\delta+d$ when $X=\{x_0,x_1,\ldots,x_m\}$}
\label{tbl:composition-products}
\begin{center}
{\scriptsize
\begin{tabular}{cclcl}
\hline\noalign{\smallskip}
\multicolumn{1}{c}{Name} & \multicolumn{1}{c}{Symbol} & \multicolumn{1}{c}{Map}
& \multicolumn{1}{c}{Operator Identity} & \multicolumn{1}{c}{Remarks} \\
\noalign{\smallskip}\svhline\noalign{\smallskip}
composition          & $c\circ d$               &
$\allseriesell\times\allseriesm\rightarrow \allseriesell$         & $F_c\circ
F_d=F_{c\circ d}$ & associative \\
modified composition & $c\modcomp d$            &
$\allseriesell\times\allseriesm\rightarrow \allseriesell$         &
$F_c\circ(I+F_d)$ & nonassociative \\
mixed composition    & $c\circ d_\delta$        &
$\allseriesell\times\allseriesdeltam\rightarrow \allseriesell$    & $F_c\circ
F_{d_\delta}=F_{c\circ d_\delta}$ & $c\circ d_\delta=c\modcomp d$ \\
group composition    & $c\circledcirc d$        & $\allseriesm
\times\allseriesm\rightarrow \allseriesm$            &
$(I+F_c)\circ(I+F_d)=I+F_{c\circledcirc d}$ & $c\circledcirc d=d+c\modcomp d$ \\
group product    & $c_\delta\circ d_\delta$ &
$\allseriesdeltam\times\allseriesdeltam\rightarrow \allseriesdeltam$  &
$F_{c_\delta}\circ F_{d_\delta}=F_{c_\delta\circ d_\delta}$ & $c_\delta\circ
d_\delta=\delta+c\circledcirc d$ \\
\noalign{\smallskip}\hline\noalign{\smallskip}
\end{tabular}
}
\end{center}
\end{table}

It is also known that the modified composition product is an ultrametric
contraction on $\allseriesm$. A summary description of the composition and
modified composition product is given in Table~\ref{tbl:composition-products}
along with other types of composition products which will be presented shortly.


\subsection{Output feedback}
\label{subsec:feedback-interconnection}

\begin{figure}[t]
\begin{center}
	\includegraphics*[scale=0.38]{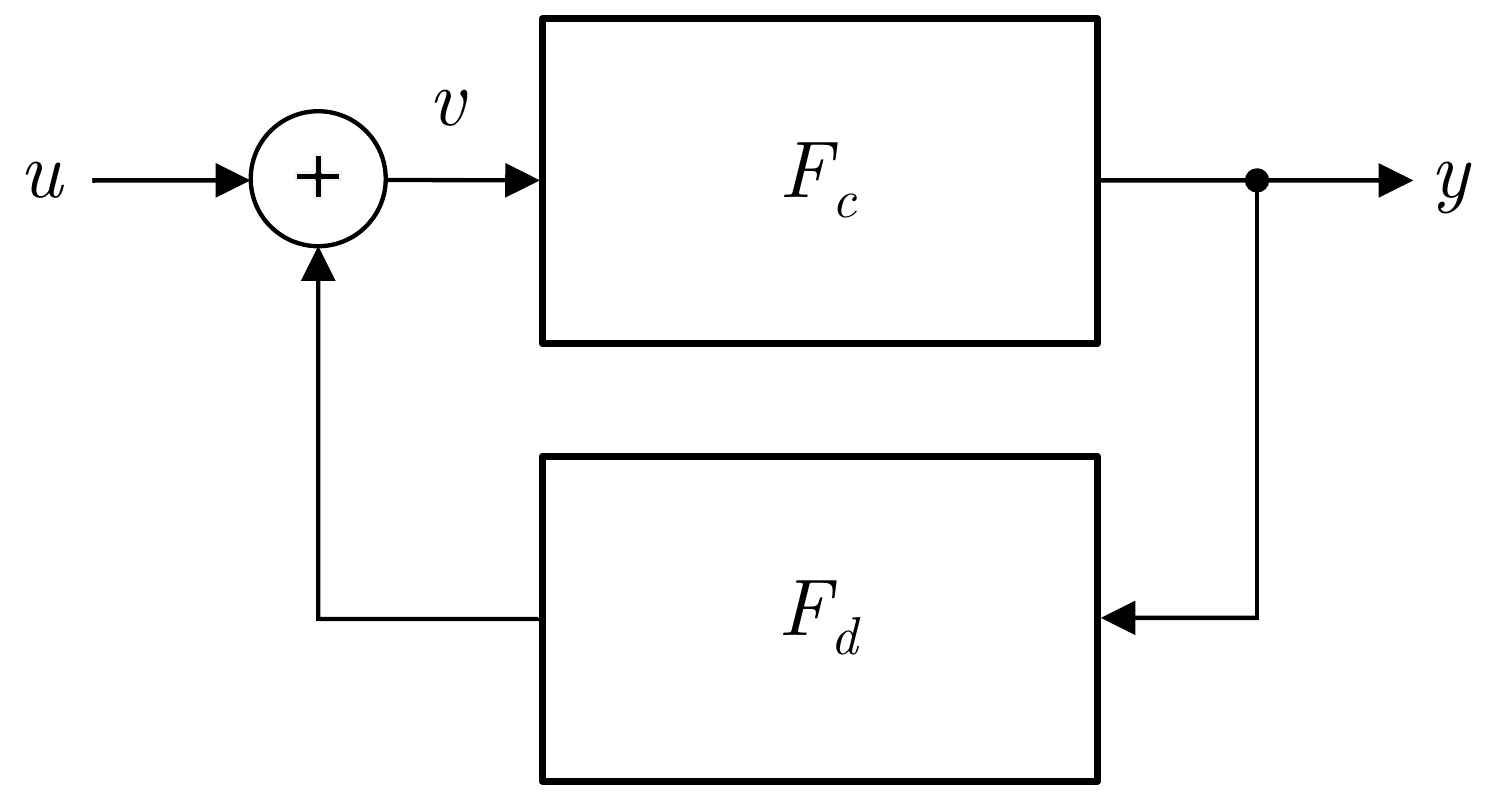}
\end{center}
\caption{Output feedback connection}
\label{fig:feedback-with-v}
\end{figure}

A central object of study in control theory is the output feedback
interconnection as shown in Figure~\ref{fig:feedback-with-v}. As with the
cascade systems discussed above, this class of interconnections is also closed
in the sense that the mapping $u\mapsto y$ always has a Fliess operator
representation. The computation of the corresponding generating series, however,
is much more involved. To see the source of the difficulty, consider the
following calculation assuming {$c,d\in\allseriesmLC$} (for the most general case,
which requires two alphabets, see \cite{Gray-etal_SCL14}).
Clearly the function $v$ in Figure~\ref{fig:feedback-with-v} must satisfy the
identity
\begdi
	v=u+F_{d\circ c}[v].
\enddi
Therefore, from $(I+F_{-d\circ c})[v] =u$ one deduces that
\begin{align*}
	v=\left(I+F_{(-d\circ c)}\right)^{-1}[u]
	  =:\left(I+F_{(-d\circ c)^{-1}}\right)[u],
\end{align*}
where $I$ is the identity operator, and the superscript ``$-1$" denotes the
composition inverse in both the operator sense and in terms of formal power
series. Thus, the generating series for the closed-loop system, denoted by the
output feedback product $c@d$, is
\begeq \label{eq:c-at-d}
	F_{c@d}[u] = F_c[v]
			  = F_{c\modcomp(-d\circ c)^{-1}}[u].
\endeq
The crux of the problem is how to compute the generating series $(-d \circ
c)^{-1}$ of the inverse operator $(I+F_{(-d\circ c)})^{-1}$. The approach
described next is based on identifying an underlying combinatorial Hopf algebra
whose antipode acts on a certain character group in such a way as to explicitly
produce this inverse generating series. This requires that
one first identifies the relevant group structures.

Consider the set of {\it unital Fliess operators}
$
	\Fliessdelta:=\{I+F_c : {c \in \allseriesmLC}\}.
$
It is convenient to introduce the symbol $\delta$ as the (fictitious) generating
series for the identity map. That is, $F_\delta:=I$ such that
$I+F_c:=F_{\delta+c}=F_{c_\delta}$ with $c_\delta:=\delta+c$. The set of all
such generating series for $\Fliessdelta$ will be denoted by {$\allseriesdeltamLC$}.
The central idea is that $(\Fliessdelta,\circ,I)$ forms a group under operator
composition
\begin{align*}
	F_{c_\delta}\circ F_{d_\delta}
		&=(I+F_c)\circ(I+F_d)=I+F_d+F_c\circ(I+F_d) \\
		&=I+F_d+F_{c\modcomp d}
		  =: F_{c_\delta\circ d_\delta},
\end{align*}
{where
$$c_\delta \circ d_\delta:=\delta+d+c\modcomp d=:\delta+c\circledcirc d.$$
(The series $c\circledcirc d$ is clearly locally convergent since all the operations employed in its
definition preserve local convergence.)}
This group will be referred to as the {\em output feedback group} of unital Fliess
operators. It is natural to think of this group as acting on an arbitrary Fliess
operator, say $F_c$, to produce another Fliess operator, as is evident in
\rref{eq:c-at-d}, by defining the right action $F_c\circ F_{d_\delta}=F_{c\circ
d_\delta}$ with $c\circ d_\delta:=c\modcomp d$. This product will be referred to
as the {\em mixed composition product} on
$\allseriesell\times\allseriesdeltam$.\footnote{The same symbol will be used for
composition on $\allseriesm$, $\allseriesdeltam$, and
$\allseriesm\times\allseriesdeltam$. It will always be clear which product is
being used since the arguments of these products have a distinct notation,
namely, $c$ versus $c_\delta$.} The following lemma gives the basic properties
of the mixed composition product.

\begle \label{le:mixed-composition-product-identities} {\rm
\cite{Gray-Duffaut-Espinosa_CDC13}}
The mixed composition product
\begin{description}
\item[(1)] is left $\re$-linear;
\item[(2)] satisfies $(c\circ d_\delta)\circ e_\delta=c\circ (d_\delta\circ
e_\delta)$ (mixed associativity);
\item[(3)] is associative.
\end{description}
\endle

\begin{proof}
$\;$\\
(1) This claim follows from the left linearity of the modified composition
product.
\newline
(2) In light of the first item it is sufficient to prove the claim only for
$c=\eta\in X^k$, $k\geq 0$.
The cases $k=0$ and $k=1$ are trivial.  Assume the claim holds up to some fixed
$k\geq 0$.
Then via item (4)  in Lemma~\ref{le:modified-product-properties} and the induction
hypothesis it follows that
\begin{align*}
	((x_0\eta)\circ d_\delta)\circ e_\delta
	&=(x_0(\eta\circ d_\delta))\circ e_\delta
		=x_0((\eta\circ d_\delta)\circ e_\delta)
		=x_0(\eta\circ (d_\delta\circ e_\delta)) \\
	&=(x_0\eta)\circ (d_\delta\circ e_\delta).
\end{align*}
In a similar fashion, apply the properties (1), (4), and (5) in
Lemma~\ref{le:modified-product-properties}  to get
\begin{align*}
	((x_i\eta)\circ d_\delta)\circ e_\delta
	&=\left[ x_i(\eta\circ d_{\delta})+x_0(d_{i}\shuffle (\eta\circ d_{\delta}))\right]\circ e_\delta \\
	&=[x_i(\eta\circ d_{\delta})]\circ e_\delta+[x_0(d_{i}\shuffle (\eta\circ d_{\delta}))]\circ e_\delta \\
	&=x_i[(\eta\circ d_{\delta})\circ e_\delta]+x_0[e_i\shuffle ((\eta\circ d_{\delta})\circ e_\delta)]+ \\
	&\hspace*{0.18in}x_0[(d_{i}\circ e_\delta)\shuffle ((\eta\circ d_{\delta})\circ e_\delta)].
\end{align*}
Now employ the induction hypothesis so that
\begin{align*}
	((x_i\eta)\circ d_\delta)\circ e_\delta
	&=x_i[\eta\circ (d_{\delta}\circ e_\delta)]+x_0[e_i\shuffle (\eta\circ (d_{\delta}\circ e_\delta))]+ \\
	&\hspace*{0.18in}x_0[(d_{i}\circ e_\delta)\shuffle (\eta\circ (d_{\delta}\circ e_\delta))] \\
	&=x_i[\eta\circ (d_{\delta}\circ e_\delta)]+x_0[(e_i+(d_{i}\circ e_\delta))\shuffle (\eta\circ (d_{\delta}\circ e_\delta))]\\
	&=x_i[\eta\circ (d_{\delta}\circ e_\delta)]+x_0[(d_\delta\circ e_\delta)_i\shuffle (\eta\circ (d_{\delta}\circ e_\delta))]\\
	&=(x_i\eta)\circ(d_\delta\circ e_\delta).
\end{align*}
Therefore, the claim holds for all $\eta\in X^\ast$, and the identity is proved.
\newline
(3) Applying item (5) in Lemma~\ref{le:modified-product-properties} and the previous
result it follows
\begin{align*}
	(c_\delta\circ d_\delta)\circ e_\delta&=\delta+e+(c\circledcirc d)\circ e_\delta\\
	&=\delta+e+(d+c\circ d_{\delta})\circ e_\delta \\
	&=\delta+e+d\circ e_{\delta}+(c\circ d_{\delta})\circ e_\delta \\
	&=\delta+d\circledcirc e+c\circ (d_{\delta}\circ e_\delta) \\
	&=(c_\delta\circ(d_\delta\circ e_\delta)).
\end{align*}
\end{proof}

Given the uniqueness of generating series of Fliess operators, {their set of generating series forms the group
$(\allseriesdeltamLC,\circ,\delta)$. In fact, it is a subgroup of the group described next since it can be shown that the composition inverse preserves local
convergence \cite{Gray-etal_SCL14}.
}
\begth \label{th:allseriesdeltam-is-group} {\rm \cite{Gray-etal_SCL14}}
The triple $(\allseriesdeltam,\circ,\delta)$ is a group.
\endth

\begin{proof}
By design, $\delta$ is the identity element of the group. Associativity of the
group product was established above.
For a fixed $c_\delta\in\allseriesdeltam$, the composition inverse,
$c_\delta^{-1}=\delta+c^{-1}$, must satisfy $c_\delta\circ c_\delta^{-1}=\delta$
and $c_\delta^{-1}\circ c_\delta=\delta$, which reduce, respectively, to
\begin{subequations}
\begin{align}
	c^{-1} &=(-c)\modcomp c^{-1} \label{eq:cdelta-right-inverse} \\
		c&=(-c^{-1})\modcomp c. \label{eq:cdelta-left-inverse}
\end{align}
\end{subequations}
It was shown in \cite{Gray-Li_05} that $e\mapsto (-c)\modcomp e$ is always a
contraction on $\allseriesm$ when viewed
as an (complete) ultrametric space and thus has a unique fixed point, $c^{-1}$.
So it follows directly that $c_\delta^{-1}$ is a right inverse of $c_\delta$,
i.e., satisfies \rref{eq:cdelta-right-inverse}. To see that this same series is
also a left inverse, first observe that \rref{eq:cdelta-right-inverse} is
equivalent to
\begeq \label{eq:left-inverse-RHS-version}
	c^{-1}\modcomp 0 +c\modcomp c^{-1}=0,
\endeq
using Lemma~\ref{le:modified-product-properties}, items (1) and (2). Substituting
\rref{eq:left-inverse-RHS-version} back into
itself where zero appears and applying \rref{eq:cmod-non-associative} gives
\begin{align*}
	c^{-1}\modcomp (c\modcomp c^{-1}+c^{-1})+c\modcomp c^{-1}&=0 \\
	(c^{-1}\modcomp c)\modcomp c^{-1}+c\modcomp c^{-1}&=0.
\end{align*}
Again from left linearity of the modified composition product it follows that
\begdi
	(c^{-1}\modcomp c+c)\modcomp c^{-1}=0.
\enddi
Finally, Lemma~\ref{le:modified-product-properties}, item (3) implies that
$c^{-1}\modcomp c+c=0$, which is equivalent to \rref{eq:cdelta-left-inverse}.
This concludes the proof.
\end{proof}

In light of the identities $c \circ \delta=c$ and
Lemma~\ref{le:mixed-composition-product-identities}, item (2), it is therefore
established that $\allseriesdeltam$ acts as a right transformation group on
$\allseriesm$. In which case, the feedback product $c@d=c\modcomp(-d\circ
c)^{-1}=c\circ (-d\circ c)_\delta^{-1}$ can be viewed as specific example of
such a right action.


\section{The Hopf Algebra of Coordinate Functions}
\label{sec:cha-crt}

In this section the Hopf algebra of coordinate functions for the group $\allseriesdeltam$ is described. This provides an explicit computational framework for computing group inverses, and thus, to calculate the feedback product as described in the previous section. The strategy is to first introduce a connected graded commutative algebra and then a compatible noncocommutative coalgebra, resulting in a connected graded commutative noncocommutative bialgebra. The connectedness property ensures then that this bialgebra is a connected graded Hopf algebra. The section ends by providing a purely inductive formula for the antipode of this Hopf algebra.


\subsection{Multivariable Hopf algebra of output feedback}
\label{ssect:multiHA}

Let $X=\{x_0,x_{1},x_{2},\ldots,x_{m}\}$ be a finite alphabet with $m+1$ letters. As usual the monoid of words is denoted by $X^*$ and includes the empty word $\mathrm{e}=\emptyset$. The degree of a word $\eta = x_{i_1} \cdots x_{i_n} \in X^*$ of length $|\eta|:=n$, where $x_{i_l} \in X$, is defined by
\begin{equation}
\label{worddegree}
	\|\eta\|:= 2|\eta|_0 + |\eta|_1.
\end{equation}
Here $|\eta|_0$ denotes the number of times the letter $x_0 \in X$ appears in $\eta$, and $|\eta|_1$ is the number letters $x_{j \neq 0} \in X$ appearing in the word $\eta$. Note that $|\mathrm{e}|=0=\|\mathrm{e}\|$.

{Recall Remark \ref{rmk:groupscheme} above. For any word $\eta \in X^\ast$ and $i=1,\ldots,m$, the {\em coordinate function} $a_\eta^i$ is defined to be the element of the dual space $\mathbb{R}^*\langle \langle X_\delta \rangle \rangle$ giving the coefficient of the $i$-th component series for the word $\eta \in X^*$, namely,
$$
	a_\eta^i(c) := ( c_i , \eta ).
$$
In this context, $a_\delta^i$ denotes the coordinate function with respect to $\delta$, where $a_\delta^i(\delta)=1$ and zero otherwise.
}
Consider the vector space $V$ generated by the coordinate functions $a^i_\eta$, where $\eta \in X^*$ and $1 \leq i \leq m$. It is turned into a polynomial algebra $H$ with unit denoted by $\mathbf{1}$. By defining the degree of elements in $H$ as $\deg(\mathbf{1}):=0$ and for $k>0$, $\eta \in X^*$
\begin{equation}
\label{degree}
	\deg(a^{k}_\eta) := 1 + \|\eta\|,
\end{equation}
$\deg(a^{k}_\eta a^{l}_\kappa) := \deg(a^{k}_\eta) + \deg(a^{l}_\kappa)$, $H$ becomes a graded connected algebra, $H:=\bigoplus_{n \ge 0} H_n$.
Note, in particular, that $\deg(a^{k}_\mathrm{e}) := 1$.

The left- and right-shift maps, ${\theta}_{j}: H \to H$ respectively $\tilde{\theta}_{j}: H \to H$ are defined by
$$
	{\theta}_{j} a^k_\eta:=a^k_{x_j\eta},
	\qquad\
	\tilde{\theta}_{j} a^k_\eta:= {a}^k_{\eta x_{j}}
$$
for $x_j \in X$, and  ${\theta}_{j}\un =\tilde{\theta}_{j}\un =0$. On products these maps act by definition as derivations
$$
	{\theta}_{j} a^k_\eta a^l_\mu:=({\theta}_{j} a^k_\eta) a^l_\mu +  a^k_\eta ({\theta}_{j}a^l_\mu),
$$
and analogously for $\tilde{\theta}_{j}$. For a word $\eta = x_{i_1} \cdots x_{i_n} \in X^*$
$$
	{\theta}_{\eta} := {\theta}_{i_1}\circ \cdots \circ {\theta}_{i_n},
	\qquad\
	\tilde{\theta}_{\eta} := \tilde{\theta}_{i_n} \circ \cdots \circ \tilde{\theta}_{i_1}.
$$
Hence, any element $a^i_\eta$, $\eta \in X^*$ can be written
$$
	a^i_\eta = {\theta}_{\eta}  a_\mathrm{e}^i = \tilde{\theta}_{\eta} a_\mathrm{e}^i.
$$

Both maps can be used to define a particular coproduct $\Delta : H \to H \otimes H$. In this approach, the right-shift map is considered first.
Later it will be shown that the left-shift map gives rise to the same coproduct. The coordinate function with respect to the empty word, $a_\mathrm{e}^l$, $1 \leq l \leq m$, is defined to be primitive
\begin{align}
\label{coproduct1a}
	\Delta a_\mathrm{e}^l := a_\mathrm{e}^l \otimes \un + \un \otimes a_\mathrm{e}^l.
\end{align}
The next step is to define $\Delta$ inductively on any $a^i_\eta$, $|\eta|>0$, by specifying intertwining relations between the map $\tilde{\theta}_{\eta} $ and the coproduct
\begin{equation}
\label{cproduct1b}
	\Delta \circ {\tilde\theta}_{i}
	:= \Big(\tilde{\theta}_{i} \otimes \id + \id \otimes \tilde{\theta}_{i}
	+  \delta_{0i} \sum_{j=1}^{m} \tilde{\theta}_{j} \otimes A_\mathrm{e}^{(j)}\Big) \circ \Delta,
\end{equation}
where $\delta_{0i}$ is the usual Kronecker delta. The map $A_\mathrm{e}^{(j)}$ for $0< j \le m$ is defined by
\begin{equation}
\label{offset}
	A_\mathrm{e}^{(j)} a^i_\eta := a^i_\eta a_\mathrm{e}^{j}.
\end{equation}
The following notation is used, $\Delta \circ {\tilde\theta}_{i}={\tilde\Theta}_{i} \circ \Delta$, where
\begin{equation}
\label{cproduct1ba}
	{\tilde\Theta}_{i} :=  \tilde{\theta}_{i} \otimes \id + \id \otimes \tilde{\theta}_{i}
	+ \delta_{0i} \sum_{j=1}^{m} \tilde{\theta}_{j} \otimes A_\mathrm{e}^{(j)},
\end{equation}
and ${\tilde\Theta}_{\eta}:=\tilde{\Theta}_{i_n} \circ \cdots \circ \tilde{\Theta}_{i_1}$ for $\eta = x_{i_1} \cdots x_{i_n} \in X^*$. The functions $a^{l}_{x_j}$ for $0<l,j \le m$ are primitive since
\begin{align*}
	\Delta a^{l}_{x_j} = \big(\tilde{\theta}_{j} \otimes \id + \id \otimes \tilde{\theta}_{j}) \circ \Delta a_\mathrm{e}^l
	   &= \big(\tilde{\theta}_{j} \otimes \id + \id \otimes \tilde{\theta}_{j})(a_\mathrm{e}^l \otimes \un
				+ \un \otimes a_\mathrm{e}^l) \\
	   &=a^{l}_{x_j} \otimes \un + \un \otimes a^{l}_{x_j},
\end{align*}
which follows from $ \tilde{\theta}_{j} \un = 0$. However, for $a^{l}_{x_0}$ the coproduct is
\allowdisplaybreaks
\begin{align}
	\Delta a^{l}_{x_0}  = {\tilde\Theta}_{0} \circ \Delta a_\mathrm{e}^{l}
	&=\Big(\tilde{\theta}_{0} \otimes \id + \id \otimes \tilde{\theta}_{0}
	+ \sum_{j=1}^{m} \tilde{\theta}_{j} \otimes A_\mathrm{e}^{(j)}\Big) \circ \Delta a_\mathrm{e}^{l}  	\nonumber \\
			&= a^{l}_{x_0} \otimes \un + \un \otimes a^{l}_{x_0}
			+  \sum_{j=1}^{m} a^{l}_{x_j} \otimes a^{j}_{\mathrm{e}}. 					\label{EXcoproduct2}
\end{align}
Observe that the coproduct is compatible with the grading. Indeed, $\deg(a^{l}_{x_0}) = 1 + 2$ and $\deg(a^{l}_{x_j} \otimes a^{j}_{\mathrm{e}})=\deg(a^{l}_{x_j}) + \deg(a^{j}_{\mathrm{e}})=1+1 + 1$ for any $j>0$. For the element $a^{l}_{x_ix_j}$, $i,j >0$, one finds the following coproduct
\allowdisplaybreaks
\begin{align}
		\Delta a^{l}_{x_ix_j} = {\tilde\Theta}_{j}\circ{\tilde\Theta}_{i} \circ \Delta a_\mathrm{e}^{l}
			&=\big(\tilde{\theta}_{j} \otimes \id + \id \otimes \tilde{\theta}_{j} \big) \circ
			\big(\tilde{\theta}_{i} \otimes \id + \id \otimes \tilde{\theta}_{i} \big)
					\circ \Delta a_\mathrm{e}^{l} 									 \nonumber \\
			&= \big(\tilde{\theta}_{j}\tilde{\theta}_{i}  \otimes \id
					+ \id \otimes  \tilde{\theta}_{j} \tilde{\theta}_{i}
					+ \tilde{\theta}_{j} \otimes  \tilde{\theta}_{i}
					+  \tilde{\theta}_{i} \otimes \tilde{\theta}_{j}\big)
				\circ \Delta a_\mathrm{e}^{l} 										 \nonumber \\
					&=a^{l}_{x_ix_j} \otimes \un + \un \otimes a^{l}_{x_ix_j}.
\end{align}
Again, this follows from $ \tilde{\theta}_{k} \un = 0$ and generalizes to any word $\eta$ in the monoid $\tilde{X}^*$ made from the reduced alphabet $\tilde{X}:=X-\{x_0\}=\{x_1,\ldots,x_m\}$. For the coproduct of $a^{l}_{x_ix_0}$, $i>0$, it follows that
\allowdisplaybreaks
\begin{align}
		\lefteqn{\Delta a^{l}_{x_ix_0} = {\tilde\Theta}_{0}\circ{\tilde\Theta}_{i} \circ \Delta a_\mathrm{e}^{l}} \nonumber \\
			&=\Big(\tilde{\theta}_{0} \otimes \id + \id \otimes \tilde{\theta}_{0}
			+ \sum_{j=1}^{m} \tilde{\theta}_{j} \otimes A_\mathrm{e}^{(j)}\Big)  \circ
			\big(\tilde{\theta}_{i} \otimes \id + \id \otimes \tilde{\theta}_{i}\big)
					\circ \Delta a_\mathrm{e}^{l} 							 		\nonumber \\
			&=\Big(\tilde{\theta}_{0} \otimes \id + \id \otimes \tilde{\theta}_{0}
			+ \sum_{j=1}^{m} \tilde{\theta}_{j} \otimes A_\mathrm{e}^{(j)}\Big)
			(a^{l}_{x_i} \otimes \un + \un \otimes a^{l}_{x_i})					 		\nonumber \\
			&=a^{l}_{x_ix_0} \otimes \un + \un \otimes a^{l}_{x_ix_0}
			+ \sum_{j=1}^{m} a^{l}_{x_ix_j} \otimes a^{j}_\mathrm{e}. 					\label{coproduct21}
\end{align}
This should be compared with the coproduct of $a^{l}_{x_0x_i}$, $i>0$
\allowdisplaybreaks
\begin{align}
		\lefteqn{\Delta a^{l}_{x_0x_i} = {\tilde\Theta}_{i}\circ{\tilde\Theta}_{0} \circ \Delta a_\mathrm{e}^{l}} \nonumber \\
			&=\big(\tilde{\theta}_{i} \otimes \id + \id \otimes \tilde{\theta}_{i}\big)  \circ
			\Big(\tilde{\theta}_{0} \otimes \id + \id \otimes \tilde{\theta}_{0}
			+ \sum_{j=1}^{m} \tilde{\theta}_{j} \otimes A_\mathrm{e}^{(j)}\Big)
					\circ \Delta a_\mathrm{e}^{l} 									\nonumber \\
			&=\big(\tilde{\theta}_{i} \otimes \id + \id \otimes \tilde{\theta}_{i}\big)
			\Big(a^{l}_{x_0} \otimes \un + \un \otimes a^{l}_{x_0}
			+ \sum_{j=1}^{m} a^{l}_{x_j} \otimes a^{j}_\mathrm{e} \Big) 						\nonumber \\
			&=a^{l}_{x_0x_i} \otimes \un + \un \otimes a^{l}_{x_0x_i}
			+ \sum_{j=1}^{m} a^{l}_{x_jx_i} \otimes a^{j}_\mathrm{e}
			+ \sum_{j=1}^{m} a^{l}_{x_j} \otimes a^{j}_{x_i}. 						\label{coproduct22}
\end{align}
Finally, the coproduct of $a^{l}_{x_0x_0}$ is calculated
\begin{align}
		\lefteqn{\Delta a^{l}_{x_0x_0} = {\tilde\Theta}_{0}\circ{\tilde\Theta}_{0} \circ \Delta a_\mathrm{e}^{l}} \nonumber \\
			&=\Big(\tilde{\theta}_{0} \otimes \id + \id \otimes \tilde{\theta}_{0}
			+ \sum_{n=1}^{m} \tilde{\theta}_{n} \otimes A_\mathrm{e}^{(n)}\Big)  \circ
			\Big(\tilde{\theta}_{0} \otimes \id + \id \otimes \tilde{\theta}_{0}
			+ \sum_{j=1}^{m} \tilde{\theta}_{j} \otimes A_\mathrm{e}^{(j)}\Big)
					\circ \Delta a_\mathrm{e}^{l}  									\nonumber \\
			&=\Big(\tilde{\theta}_{0} \otimes \id + \id \otimes \tilde{\theta}_{0}
			+ \sum_{n=1}^{m} \tilde{\theta}_{n} \otimes A_\mathrm{e}^{(n)}\Big)
			\Big(a^{l}_{x_0} \otimes \un + \un \otimes a^{l}_{x_0}
			+ \sum_{j=1}^{m} a^{l}_{x_j} \otimes a^{j}_\mathrm{e} \Big) 						\nonumber \\
			&=a^{l}_{x_0x_0} \otimes \un + \un \otimes a^{l}_{x_0x_0}
			+ \sum_{j=1}^{m} a^{l}_{x_jx_0} \otimes a^{j}_\mathrm{e}
			+ \sum_{n=1}^{m} a^{l}_{x_0x_n} \otimes a^{n}_\mathrm{e} 					\nonumber \\
			& \qquad + \sum_{j=1}^{m} a^{l}_{x_j} \otimes a^{j}_{x_0}
			+ \sum_{n,j=1}^{m} a^{l}_{x_jx_n} \otimes a^{j}_\mathrm{e}a^{n}_\mathrm{e}. 	\label{coproduct23}
\end{align}

The coproduct $\Delta$ is extended multiplicatively to all of $H$, and the unit $\un$ is defined to be $\Delta(\un):=\un \otimes \un$.
{In particular, $\Delta a^{l}_{\eta}\in V \otimes H$ as opposed to simply $H\otimes H$ due to the left linearity of the mixed composition product as shown in
Lemma~\ref{le:mixed-composition-product-identities}. This has interesting practical consequences for the antipode calculation as described in \cite{Duffaut-etal_submit2014}.
Namely, there is a significant difference between the computational efficiencies of the two antipode recursions in \eqref{antipode2}.}

\begin{theorem}\label{thm:MIMO-HA}
The algebra $H$ with the multiplicatively extended coproduct
\begin{equation}
\label{coproductMIMO}
	\Delta a^{l}_{\eta} := {\tilde\Theta}_{\eta}( a_\mathrm{e}^{l} \otimes \un + \un \otimes a_\mathrm{e}^{l}) \in V \otimes H
\end{equation}
is a connected graded commutative noncocommutative Hopf algebra.
\end{theorem}

\begin{proof}
Since $H$ is connected, graded and commutative by construction, the antipode can be calculated recursively via one of the recursions in \eqref{antipode2}. It is clear as well that the coproduct \eqref{coproductMIMO} is noncocommutative. Hence, coassociativity remains to be checked. This is done inductively with respect to the length of the word $\eta \in X^*$.
First observe that
$$
	\Delta(a_{\eta x_i}^k) = \Delta \circ \tilde\theta_{i} (a_{\eta}^k)
					 = \tilde\Theta_{i} \circ  \Delta(a_{\eta}^k).
$$
This implies that
\begin{align}
	\lefteqn{(\Delta \otimes \id) \circ \Delta(a_{\eta x_i}^k)
	= (\Delta \otimes \id)  \circ \tilde \Theta_{i} \circ \Delta(a_{\eta}^k)}					\nonumber \\
	&= \Big( \Delta  \circ \tilde{\theta}_{i} \otimes \id +  (\id \otimes \id) \circ \Delta\otimes \tilde{\theta}_{i}
	+ \delta_{0i}\sum_{j=1}^{m} \Delta \circ  \tilde{\theta}_{j} \otimes A_\mathrm{e}^{(j)}\Big) \circ \Delta(a_{\eta}^k)\\
	&=\Big(\tilde{\Theta}_{i} \otimes \id +  \id \otimes \id \otimes \tilde{\theta}_{i}
	+ \delta_{0i} \sum_{j=1}^{m} \tilde{\Theta}_{j} \otimes A_\mathrm{e}^{(j)}\Big)
	(\Delta \otimes \id)  \circ \Delta(a_{\eta}^k)										\nonumber \\
	&=\Big(\tilde{\theta}_{i} \otimes \id \otimes \id
		+ \id \otimes \tilde{\theta}_{i} \otimes \id
		+  \id \otimes \id \otimes \tilde{\theta}_{i} 									\nonumber \\
	& \qquad + \delta_{0i}\sum_{j=1}^{m} \tilde{\theta}_{j} \otimes A_\mathrm{e}^{(j)}\otimes \id
	 	+ \delta_{0i}\sum_{j=1}^{m} \tilde{\theta}_{j} \otimes \id \otimes A_\mathrm{e}^{(j)} 	\nonumber \\
	& \qquad + \delta_{0i}\sum_{j=1}^{m}  \id \otimes \tilde{\theta}_{j} \otimes A_\mathrm{e}^{(j)}\Big)
		(\id \otimes \Delta)  \circ \Delta(a_{\eta}^k)									\nonumber \\
	&= \Big(\tilde{\theta}_{i} \otimes \id \otimes \id + \id \otimes \tilde{\Theta}_{i}
		+ \delta_{0i}\sum_{j=1}^{m} \tilde{\theta}_{j} \otimes \big(A_\mathrm{e}^{(j)}\otimes \id
		+ \id \otimes A_\mathrm{e}^{(j)}\big)\Big)(\id \otimes \Delta)  \circ \Delta(a_{\eta}^k)	\nonumber \\
	&= (\id \otimes \Delta)\circ
	\Big(\tilde{\theta}_{i} \otimes \id + \id \otimes \tilde{\theta}_{i}
		+  \delta_{0i}\sum_{j=1}^{m} \tilde{\theta}_{j} \otimes A_\mathrm{e}^{(j)}\Big) \circ \Delta(a_{\eta}^k) \nonumber \\
	&= (\id \otimes \Delta)\circ\Delta(a_{\eta x_i}^k).
\end{align}
The identity
$$
	\Delta \circ A_\mathrm{e}^{(i)}
	=
	\big(A_\mathrm{e}^{(i)}\otimes \id  + \id \otimes A_\mathrm{e}^{(i)}\big)\circ\Delta
$$
was used above, which follows from $A_\mathrm{e}^{(l)} a_{\eta}^k = a_\mathrm{e}^l a_{\eta}^k$ and $a_\mathrm{e}^l$ being primitive for all $0 < l \le m$.
\end{proof}

\begin{remark}\label{simplecoprod}
Note that the coproduct \eqref{coproductMIMO} can be simplified
\begin{equation}
\label{coproductMIMOa}
	\Delta a^{l}_{\eta} := {\tilde\Theta}_{\eta}( a_\mathrm{e}^{l} \otimes \un) + \un \otimes a^{l}_{\eta},
\end{equation}
which follows  from $ \tilde{\theta}_{k} \un = 0$ and the form of $\tilde\Theta_{i}$.
\end{remark}
A variant of Sweedler's notation is used for the reduced coproduct, i.e., $\Delta'(a^{l}_{\eta})= \sideset{}{'}\sum a^{l}_{\eta'} \otimes a^{l'}_{\eta''}$, as well as for the full coproduct
$$
	\Delta(a^{l}_{\eta}) = \sum a^{l}_{\eta_{(1)}} \otimes a^{l'}_{\eta_{(2)}}
				      = a^{l}_{\eta} \otimes {\bf{1}} + {\bf{1}} \otimes a^{l}_{\eta} + \Delta'(a^{l}_{\eta}).
$$
Connectedness of $H$ implies that its antipode $S: H \to H$ can be calculated using \eqref{antipode2}, namely,
\begin{equation}
\label{antipoderecursions}
	S a^{l}_{\eta} = -a^{l}_{\eta} - \sideset{}{'}\sum S(a^{l}_{\eta'})a^{l'}_{\eta''}
			= -a^{l}_{\eta} - \sideset{}{'}\sum a^{l}_{\eta'}S(a^{l'}_{\eta''}).
\end{equation}
A few examples are given next. The coproduct \eqref{coproduct1a} implies for the elements $a^{k}_\mathrm{e}$ that $S a^{k}_\mathrm{e} = - a^{k}_{\mathrm{e}}$. For $0< j,k,l \le m$
\begin{align}	
	S a^{k}_{x_j} = - a^{k}_{x_j},
	\qquad							
	S a^{l}_{x_0} = - a^{l}_{x_0}
					+ \sum_{i=1}^m a^{l}_{x_i} a^{i}_{\mathrm{e}}. 	
	\label{antipodeexample}
\end{align}

The next theorem uses the coproduct formula \eqref{coproductMIMO} to provide an alternative formula for the antipode of $H$.

\begin{theorem}\label{thm:antipode}
For any nonempty word $\eta= x_{i_1} \cdots x_{i_l}$, the antipode $S: H \to H$ satisfies
\begin{equation}
\label{antipodeSpecial}
	S a^k_\eta =  (-1)^{|\eta|+1}\tilde\Theta'_{\eta}(a^k_\mathrm{e}),
\end{equation}
where
\begin{equation}
\label{A1}
	 \tilde\Theta'_{\eta}:= \tilde\theta'_{i_l} \circ \cdots \circ \tilde\theta'_{i_1}
\end{equation}
and
\begin{equation}
\label{A2}
	\tilde\theta'_{l}:= - \tilde \theta_{l}
				+ \delta_{0l}\sum_{j=1}^{m} a_\mathrm{e}^{j} \tilde{\theta}_{j}.
\end{equation}
\end{theorem}

\begin{proof}
The claim is equivalent to saying that
$$
	S \circ \tilde{\theta}_{\eta} = - \tilde \theta'_{i_l} \circ S \circ \tilde \theta_{i_{l-1}} \circ \cdots \circ  \tilde\theta_{i_1}
$$
for the word $\eta=x_{i_1} \cdots x_{i_l} \in X^*$. The proof is via induction on the degree of $a^k_\eta$. The degree one case is excluded by assumption as it corresponds to $S a^{k}_\mathrm{e} = - a^{k}_{\mathrm{e}}$. For degree two, three, four and five it is quickly verified for $i>0$ that
$$
	 Sa^k_{x_i} = S \circ  \tilde\theta_{i} a^k_\mathrm{e}
	 = - \tilde\theta'_{i} \circ S a^k_\mathrm{e} = \tilde\theta'_{i} a^k_\mathrm{e}=- \tilde \theta_{i}a^k_\mathrm{e}= -a^k_{x_i}
$$
and
$$
	Sa^k_{x_0} = S \circ  \tilde\theta_{0} a^k_\mathrm{e}
	=\tilde\theta'_{0}a^k_\mathrm{e}
	=-a^k_{x_0}+\sum_{i=1}^{m}a^k_{x_i}a^i_\mathrm{e},
$$
which coincide with \eqref{antipodeexample}. For $j>0$
\begin{align*}
	Sa^k_{x_jx_0}   = S \circ \tilde\theta_{0}\circ\tilde\theta_{j}a^k_\mathrm{e}
			       &= (-1)\tilde\theta'_{0} \circ\tilde\theta'_{j}(a^k_\mathrm{e})
			       =-a^k_{x_jx_0}+\sum_{i=1}^{m}a^k_{x_jx_i}a^i_\mathrm{e}\\
			       &= \tilde\theta'_{0} \circ\tilde\theta'_{j}S(a^k_\mathrm{e})\\
			       &= -\tilde\theta'_{0} \circ S \circ \tilde\theta_{j}(a^k_\mathrm{e}).
\end{align*}
For degree five
\begin{align*}
	\lefteqn{Sa^k_{x_0x_0} = S\circ\tilde\theta_{0}\circ\tilde\theta_{0}a^k_\mathrm{e}
						= - \tilde\theta'_{0}\circ\tilde\theta'_{0}(a^k_\mathrm{e})}\\
				&= - \Big(-\tilde \theta_{0} + \sum_{j=1}^{m} a_\mathrm{e}^{j} \tilde{\theta}_{j}\Big)
					\Big(-\tilde \theta_{0} + \sum_{n=1}^{m} a_\mathrm{e}^{n} \tilde{\theta}_{n}\Big) (a^k_\mathrm{e})\\
				&=\Big( -\tilde \theta_{0}\tilde \theta_{0}
					+ \sum_{n=1}^{m} a_{x_0}^{n} \tilde{\theta}_{n}
					+ \sum_{n=1}^{m} a_\mathrm{e}^{n} \tilde \theta_{0} \tilde{\theta}_{n}
				 	+  \sum_{j=1}^{m} a_\mathrm{e}^{j} \tilde{\theta}_{j} \tilde \theta_{0}
				 	- \sum_{j=1}^{m}\sum_{n=1}^{m}  a_\mathrm{e}^{j} \tilde{\theta}_{j} a_\mathrm{e}^{n} \tilde{\theta}_{n}\Big)(a^k_\mathrm{e})\\
				 &=\Big( - \tilde \theta_{0}\tilde \theta_{0}
				 + \sum_{n=1}^{m} a_{x_0}^{n} \tilde{\theta}_{n}
				 + \sum_{n=1}^{m} a_\mathrm{e}^{n} \tilde \theta_{0} \tilde{\theta}_{n}
				 +  \sum_{j=1}^{m} a_\mathrm{e}^{j} \tilde{\theta}_{j} \tilde \theta_{0} \\
				 & \qquad\qquad\
				 - \sum_{j=1}^{m}\sum_{n=1}^{m}  a_\mathrm{e}^{j} a_{x_j}^{n} \tilde{\theta}_{n}
				 - \sum_{j=1}^{m}\sum_{n=1}^{m}  a_\mathrm{e}^{j}  a_\mathrm{e}^{n} \tilde{\theta}_{j}\tilde{\theta}_{n}\Big)(a^k_\mathrm{e})\\
				&= -a^k_{x_0x_0}
				+ \sum_{n=1}^{m} a^k_{x_n}a_{x_0}^{n}
				+\sum_{n=1}^{m}a^k_{x_nx_0}a^n_\mathrm{e}	
					 +\sum_{n=1}^{m}a^k_{x_0x_n}a^n_\mathrm{e}\\
				& \qquad\qquad\  	
					 - \sum_{j=1}^{m}\sum_{n=1}^{m}  a^k_{x_n}a_{x_j}^{n} a_\mathrm{e}^{j}
					  - \sum_{j=1}^{m}\sum_{n=1}^{m}   a_{x_nx_j}^{n} a_\mathrm{e}^{j}  a_\mathrm{e}^{n}.
\end{align*}

Recall that Sweedler's notation for the reduced coproduct is in use. It is assumed that the theorem holds up to degree $n \geq 2$. Recall that $\tilde \theta'_{i_l}$ are derivations on $H$ and that for the augmentation ideal projector $P$ it holds that $P \un=0$. Working with the second recursion in \eqref{antipoderecursions} one finds for $\deg(a_\eta^k)=n+1$, $\eta=x_{i_1} \cdots x_{i_l}=\bar{\eta}x_{i_l} \in X^*$ that
\begin{align*}
	S a_\eta^k &= m_H \circ (P \otimes S) \circ \Delta a_\eta^k\\
			 &= m_H \circ \big(P \circ \tilde\theta_{i_l} \otimes S + P \otimes S \circ \tilde\theta_{i_l}
			 + \delta_{0i_l} \sum_{n=1}^m P \circ \tilde\theta_{n} \otimes S \circ A^{(n)}_\mathrm{e}\big) \circ \Delta a_{\bar{\eta}}^k\\
			 &= m_H \circ \big(P \circ \tilde\theta_{i_l} \otimes S\big) \circ (a_{\bar{\eta}}^k \otimes \un + \Delta' a_{\bar{\eta}}^k)
			 + m_H \circ \big(P \otimes S \circ \tilde\theta_{i_l}\big)  \circ \Delta' a_{\bar{\eta}}^k\\
			 & \qquad\qquad\  	
			 +m_H \circ \big( \delta_{0i_l} \sum_{n=1}^m P \circ\tilde\theta_{n} \otimes S \circ A^{(n)}_\mathrm{e}\big)
			 \circ (a_{\bar{\eta}}^k \otimes \un + \Delta' a_{\bar{\eta}}^k).
\end{align*}
The critical term is
$$
	m_H \circ (P \otimes S \circ \tilde\theta_{i_l})  \circ \Delta' a_{\bar{\eta}}^k
	=m_H \circ (P \otimes S \circ \tilde\theta_{i_l})  \circ \sideset{}{'}\sum a^{l}_{\eta'} \otimes a^{l'}_{\eta''}.
$$
Since $\deg (\tilde\theta_{i_l} a^{l'}_{\eta''}) < n+1$, it can be written as
$$
	m_H \circ (P \otimes S \circ \tilde\theta_{i_l})  \circ \Delta' a_{\bar{\eta}}^k
	=-m_H \circ (P \otimes \tilde\theta'_{i_l} \circ S)  \circ \Delta' a_{\bar{\eta}}^k.
$$
This yields
\begin{align*}
	\lefteqn{S a_\eta^k = m_H \circ (P \otimes S) \circ \Delta a_\eta^k}\\
			 &= m_H \circ \big(P \circ \tilde\theta_{i_l} \otimes S -  P \otimes \tilde\theta'_{i_l} \circ S
			 + \delta_{0i_l} \sum_{n=1}^m P \circ \tilde\theta_{n} \otimes S \circ  A^{(n)}_\mathrm{e}\big) \circ \Delta a_{\bar{\eta}}^k\\
			 &= m_H \circ \big(\tilde\theta_{i_l} \otimes \id
			 			+ \id \otimes \tilde\theta_{i_l}
						-  \id \otimes  \delta_{0i_l} \sum_{n=1}^m A^{(n)}_\mathrm{e}\tilde\theta_{n}
			 			- \delta_{0i_l} \sum_{n=1}^m  \tilde\theta_{n} \otimes   A^{(n)}_\mathrm{e}\big) \circ (P \otimes S)\circ \Delta a_{\bar{\eta}}^k\\
			&= - \big(-\tilde\theta_{i_l} + \delta_{0i_l} \sum_{n=1}^m A^{(n)}_\mathrm{e}\tilde\theta_{n} \big)
			m_H \circ (P \otimes S) \circ \Delta a_{\bar\eta}^k\\
			&= - \tilde\theta'_{i_l}Sa_{\bar\eta}^k,
\end{align*}
which proves the theorem.
Note that the next to the last equality used the fact that the $\tilde\theta_{i}$ are derivations on $H$.
\end{proof}


{{
\begin{remark} Consider the case where $m=1$ in Theorem~\ref{thm:MIMO-HA}. That is, the alphabet $X:=\{x_0,x_1\}$, and the Hopf algebra $H$ is generated by the coordinate functions $a_{\eta}$, $\eta \in X^*$. Note that the upper index on the coordinate functions can be dismissed as $m=1$. The element $a_{\eta} \in H$ has the coproduct defined in terms of $\Delta \circ {\tilde\theta}_{i}={\tilde\Theta}_{i} \circ \Delta$, $i=0,1$, where ${\tilde\Theta}_{1}=\tilde{\theta}_{1} \otimes \id + \id \otimes \tilde{\theta}_{1}$ and ${\tilde\Theta}_{0} :=  \tilde{\theta}_{0} \otimes \id + \id \otimes \tilde{\theta}_{0} + \tilde{\theta}_1 \otimes A_\mathrm{e}.$ The antipode $S: H \to H$ for any nonempty word $\eta= x_{i_1} \cdots x_{i_l}$ is given by $S a_\eta =  (-1)^{|\eta|+1}\tilde\Theta'_{\eta}(a_\mathrm{e}),$ where $\tilde\Theta'_{\eta}:= \tilde\theta'_{i_p} \circ \cdots \circ \tilde\theta'_{i_1}$ and $\tilde\theta'_{l}:= - \tilde \theta_{l} + \delta_{0l} a_\mathrm{e} \tilde{\theta}_{1}$, $l=0,1$. Here $|\eta |:= 2|\eta|_0 + |\eta|_1$, where $|\eta|_0$ denotes the number of times the letter $x_0$ appears in $\eta$, and $|\eta|_1$ is the number of times the letter $x_{1}$ is appearing in the word $\eta$. One can verify directly that this Hopf algebra coincides with the single-input/single-output (SISO) feedback Hopf algebra described in \cite{Gray-Duffaut_Espinosa_SCL11}. The reader is also referred to \cite{Duffaut-etal_submit2014,Ebrahimi-Fard-Gray_IMRN17} for more details. This connection between Hopf algebras will be studied further in future work regarding the multivariable (MIMO) case as described in \cite{Gray-etal_SCL14}.
\end{remark}}}

Finally, returning to the antipode recursions in \eqref{antipoderecursions} one realizes quickly the intricacies that result from the signs of the different terms. Surprisingly, the computational aspects of the two formulas are rather different. It turns out that the rightmost recursion is optimal in the sense that its expansion is free of cancellations \cite{Duffaut-etal_submit2014}. This triggers immediately the question whether the antipode formula \eqref{antipodeSpecial} shares similar properties, which is answered by the next result.

\begin{proposition}
The antipode formula \eqref{antipodeSpecial} is free of cancellations.
\end{proposition}

\begin{proof}
First recall that the algebra of coordinate functions is polynomially free by construction. Then the absence of cancellations follows from looking at \eqref{A1} and \eqref{A2} and noting that $\tilde\theta_{i}\tilde\theta_{j}\neq \tilde\theta_{j}\tilde\theta_{i}$.

\end{proof}


\section{Towards Discretization}
\label{sec:towards-discretization}

This section lays the foundation for a discrete-time analogue of the continuous-time Fliess operator theory described in the previous sections.
{The starting point is the introduction of a {\em discrete-time} Fliess operator, where the basic idea is to replace the iterated integrals in
\rref{eq:Fliess-operator-defined} with iterated sums. This concept was originally exploited in \cite{Gray-etal_ACC16,Gray-etal_NM} to provide numerical
approximations of continuous-time Fliess operators. The main results were developed without any a priori assumption regarding the existence of a state space realization.
It was shown, however, that discrete-time Fliess operators are realizable by the class of state space realizations which are rational in the input and affine in the state whenever the generating series is rational. Some specific examples of this will be given here.}

\begin{figure}[t]
\begin{center}
\includegraphics[scale=0.3]{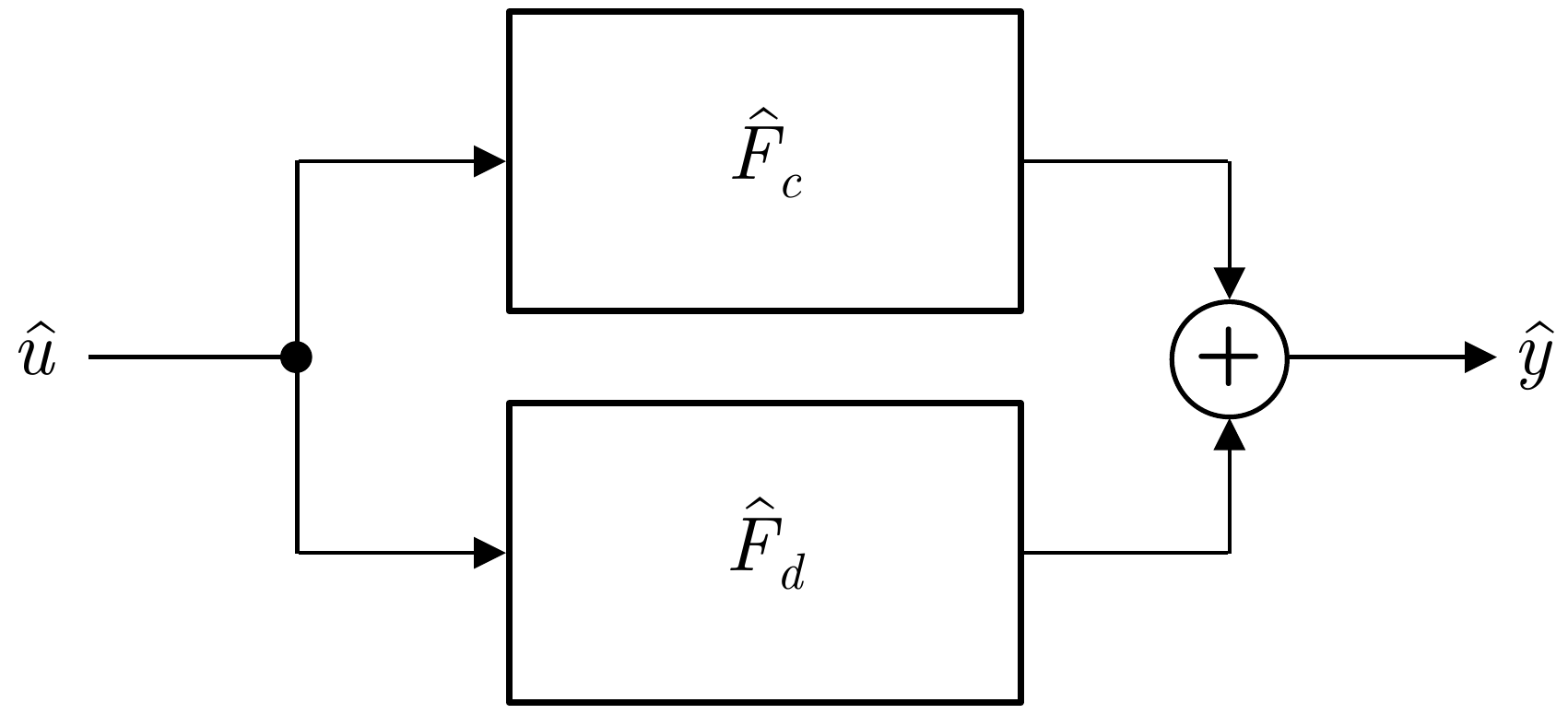}
\;\;\;\;\includegraphics[scale=0.3]{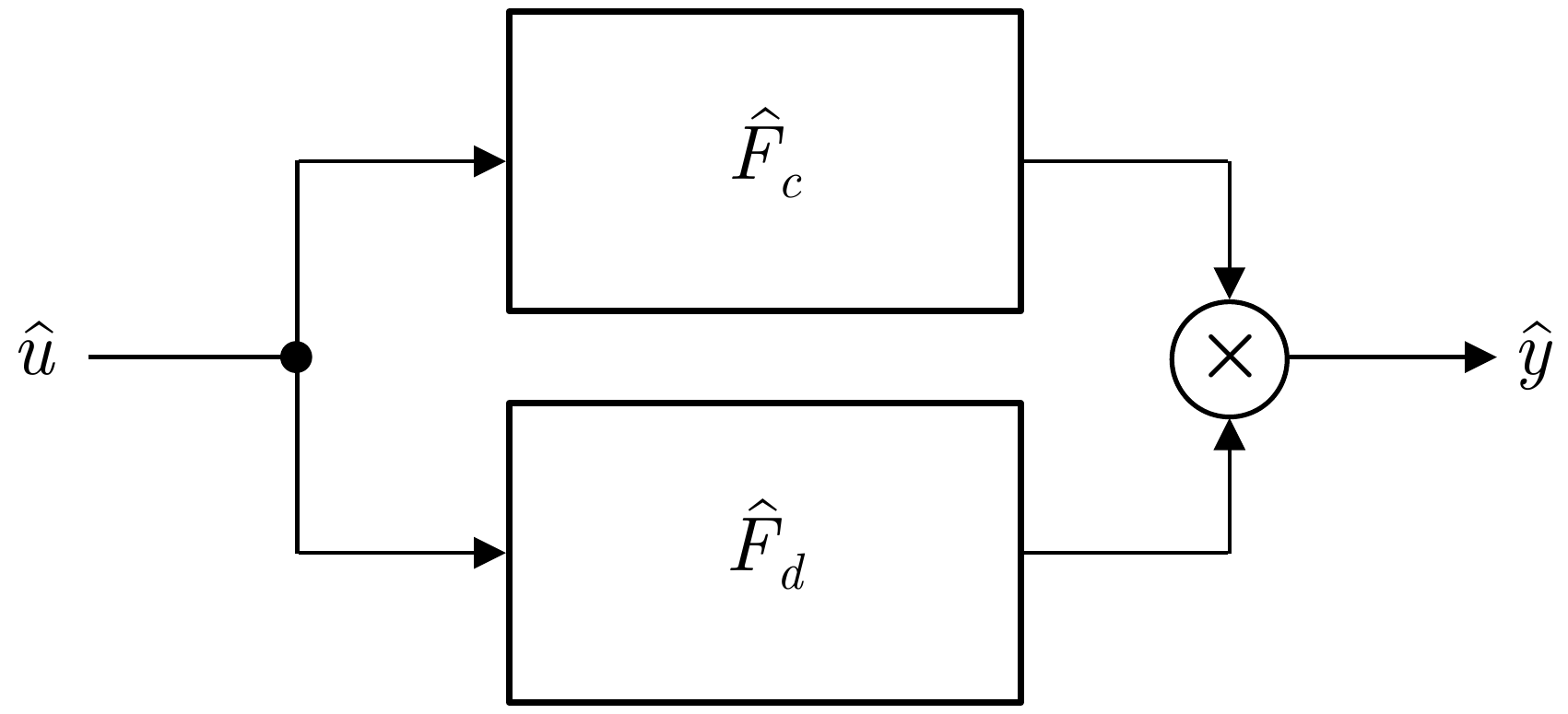}
\caption{Parallel sum (left) and parallel product (right) interconnections of two discrete-time Fliess operators.}
\label{fig:parallel}
\end{center}
\end{figure}

The main focus of this section is on parallel interconnections of discrete-time Fliess operators as shown in Figure~\ref{fig:parallel}.
{In the continuous-time theory presented earlier, it was evident that virtually all the results about interconnections flow from
the shuffle algebra, which is induced by the parallel product interconnection. The hypothesis here is that an analogous situation holds in the discrete-time case.
So it will be shown that the parallel product of discrete-time Fliess operators induces a {\em quasi}-shuffle algebra on the set of generating series. Given the natural suitability of rational generating series for discrete-time realization theory, a natural question to pursue is whether rationality is preserved under the quasi-shuffle product.
The question was affirmatively answered in \cite{Houseaux-etal_03} but without proof. So here a complete proof
will be given.}


\subsection{Discrete-time Fliess operators}
\label{subsec:discete-time-Fliess-operators}

The set of admissible inputs for discrete-time Fliess operators will be drawn from the real sequence space
\begdi
	l_\infty^{m+1}[N_0]:=\{\hat{u}=(\hat{u}(N_0),\hat{u}(N_0+1),\ldots):{ \exists \hat{R}_u \mbox{ with }0\leq}\abs{\hat{u}(N)}<\hat{R}_{\hat{u}}<\infty,\, \forall N\geq N_0 \},
\enddi
where $\hat{u}=[\hat{u}_0,\hat{u}_1,\ldots,\hat{u}_m]^T$ and $\abs{\hat{u}(N)}:=\max_{i=0,1,\ldots,m}\abs{\hat{u}_i(N)}$. In which case, $\norm{\hat{u}}_\infty:=$ \linebreak $\sup_{N\geq N_0}\abs{\hat{u}(N)}$ is always finite. Define a ball of radius $\hat{R}$ in $l_\infty^{m+1}[N_0]$ as
\begdi
	B_\infty^{m+1}[N_0](\hat{R})=\{\hat{u}\in l_\infty^{m+1}[N_0]: \norm{\hat{u}}_\infty\leq \hat{R} \}.
\enddi
The subset of finite sequences over $[N_0,N_f]$ is denoted by $B_\infty^{m+1}[N_0,N_f](\hat{R})$.
{That is, $\hat{u}\in B_\infty^{m+1}[N_0,N_f](\hat{R})$ if $\max_{N\in[N_0,N_f]} \abs{\hat{u}(N)}\leq \hat{R}$.}
The following definition is of central importance.

\begde \cite{Gray-etal_ACC16,Gray-etal_NM}
For any $c\in\allseriesell$, the corresponding {\bf discrete-time Fliess
operator} is
\begeq \label{eq:DT-Fliess-operator-defined}
	\hat{y}(N)=\hat{F}_c[\hat{u}](N)=\sum_{\eta\in X^\ast}(c,\eta)S_{\eta}[\hat{u}](N),
\endeq
where $\hat{u}\in l_\infty^{m+1}[1]$, $N\geq 1$, and the iterated sum for any $x_i\in X$ and
$\eta\in X^\ast$ is defined inductively by
\begeq \label{eq:interated-sum-approx-for-iterated-integral}
	S_{x_i\eta}[\hat{u}](N)=\sum_{k=1}^N \hat{u}_{i}(k)S_{\eta}[\hat{u}](k)
\endeq
with $S_\emptyset[\hat{u}](N):=1$.
\endde

The following lemma will be used for providing sufficient conditions for the convergence of such operators.

\begle \label{le:Seta-bound} \cite{Gray-etal_ACC16,Gray-etal_NM}
If $\hat{u}\in B_\infty^{m+1}[1](\hat{R})$ then for any $\eta\in X^\ast$ and
$N\geq 1$
$$
	\abs{S_\eta[\hat{u}](N)}\leq \hat{R}^{\abs{\eta}} {N-1+\abs{\eta} \choose \abs{\eta}} \leq 2^{N-1}(2\hat{R})^{\abs{\eta}}. \\[0.1in]
$$
\endle

\begin{proof}
If $\eta=x_{i_j}\cdots x_{i_1} \in X^*$ then for any $N\geq 1$
\allowdisplaybreaks
\begin{align*}
	\abs{S_\eta[\hat{u}](N)}
	&= \abs{\sum_{k_j=1}^N \hat{u}_{i_j}(k_j)\sum_{k_{j-1}=1}^{k_j}
	\hat{u}_{i_{j-1}}(k_{j-1})\cdots \sum_{k_1=1}^{k_2} \hat{u}_{i_1}(k_1)} \\
	&\leq \sum_{k_j=1}^N \abs{\hat{u}_{i_j}(k_j)}\sum_{k_{j-1}=1}^{k_j}
	|\hat{u}_{i_{j-1}}(k_{j-1})|\cdots \sum_{k_1=1}^{k_2}
	\abs{\hat{u}_{i_1}(k_1)} \\
	&\leq \hat{R}^{\abs{\eta}} \sum_{k_j=1}^N \sum_{k_{j-1}=1}^{k_j} \cdots
	\sum_{k_1=1}^{k_2} 1 = \hat{R}^{\abs{\eta}} {N-1+\abs{\eta} \choose \abs{\eta}},
\end{align*}
using the fact that the final nested sum above has ${N-1+\abs{\eta} \choose \abs{\eta}}$ terms \cite{Butler-Karasik_10}. The remaining inequality is standard.
\end{proof}

Since the upper bound on $\abs{S_\eta[\hat{u}](N)}$ in this lemma is achievable, it is not difficult to see that when the generating series $c$ satisfies the growth bound \rref{eq:local-convergence-growth-bound}, the series \rref{eq:DT-Fliess-operator-defined} defining $\hat{F}_c$ can diverge. For example, if $(c,\eta)=K_cM_c^{\abs{\eta}}\abs{\eta}!$ for all $\eta\in X^\ast$, and $\hat{u}_i(N)=\hat{R}$, $N\geq 1$, $i=0,1,\ldots,m$ then
\allowdisplaybreaks
\begin{align*}
	F[\hat{u}](N)
	&= \sum_{\eta\in X^\ast}K_cM_c^{\abs{\eta}}\abs{\eta}!\, \hat{R}^{\abs{\eta}}{N-1+\abs{\eta} \choose \abs{\eta}} \\
	&=K_c\sum_{j=0}^\infty (M_c(m+1)\hat{R})^j j! {N-1+j \choose j}.
\end{align*}
Since $\lim_{j\rightarrow \infty}{N-1+j \choose j}=1$, this series diverges even when $\hat{R}<1/M_c(m+1)$. The next theorem shows that this
problem is averted when $c$ satisfies the stronger growth condition \rref{eq:global-convergence-growth-bound}.

\begth \cite{Gray-etal_ACC16,Gray-etal_NM}
Suppose $c\in\allseriesell$ has coefficients which satisfy
\begdi
{\abs{(c,\eta)}\le K_c M_c^{|\eta|},\quad \forall\eta\in X^{\ast}.}
\enddi
Then there exists a real number $\hat{R}>0$ such that for each $\hat{u}\in B_\infty^{m+1}[1](\hat{R})$, the series \rref{eq:DT-Fliess-operator-defined} converges absolutely for any $N\geq 1$.
\endth

\begin{proof}
Fix $N\geq 1$. From the assumed coefficient bound and Lemma~\ref{le:Seta-bound}, it follows that
\allowdisplaybreaks
\begin{align*}
	\abs{\hat{F}_c(\hat{u})(N)}&\leq \sum_{j=0}^\infty \sum_{\eta\in X^j}
	\abs{(c,\eta)} \abs{S_\eta[\hat{u}](N)}
	\leq\sum_{j=0}^\infty K_c(M_c(m+1))^j\, 2^{N-1}(2\hat{R})^j \\
	&=\frac{K_c2^{N-1}}{1-2M_c(m+1)\hat{R}},
\end{align*}
provided $\hat{R}<1/2M_c(m+1)$.
\end{proof}

The final convergence theorem shows that the restriction on the norm of $\hat{u}$ can be removed if an even more stringent growth condition is imposed on $c$.

\begth \cite{Gray-etal_ACC16,Gray-etal_NM}
Suppose $c\in\allseriesell$ has coefficients which satisfy
\begdi \label{eq:exp-growth-bound}
	\abs{(c,\eta)}\leq K_cM_c^{\abs{\eta}}\frac{1}{\abs{\eta}!},\quad \forall\eta\in X^\ast
\enddi
for some real numbers $K_c,M_c>0$. Then for every $\hat{u}\in l_\infty^{m+1}[1]$, the series \rref{eq:DT-Fliess-operator-defined} converges absolutely for any $N\geq 1$.
\endth

\begin{proof}
Following the same argument as in the proof of the previous theorem, it is clear for any $\hat{u}\in l_\infty^{m+1}[1]$ and $N\geq 1$ that
\begdi
	\abs{\hat{F}_c(\hat{u})(N)}
	\leq \sum_{j=0}^\infty K_c(M_c(m+1))^j \frac{1}{j!}\,2^{N-1}(2\|\hat{u}\|_\infty)^j
	=K_c2^{N-1}\expup^{2M_c(m+1)\|\hat{u}\|_\infty}.
\enddi
\end{proof}

Assuming the analogous definitions for local convergence (LC) and global convergence (GC) of the operator $\hat{F}_c$, note the incongruence between the convergence conditions for continuous-time and discrete-time Fliess operators as summarized in Table~\ref{tbl:convergence-cond-summary}. In each case, for a fixed series $c$, the sense in which the discrete-time Fliess operator $\hat{F}_c$ converges is {\em weaker} than that for $F_c$. The source of this dichotomy is the observation in Lemma~\ref{le:Seta-bound} that iterated sums of $\hat{u}$ do not grow as a function of word length like $\hat{R}^{\abs{\eta}}/\abs{\eta}!$, which is the case for iterated integrals, but at a potentially faster rate. On the other hand, it is well known that rational series have coefficients that grow as in \rref{eq:global-convergence-growth-bound}. Thus, as indicated in Table~\ref{tbl:convergence-cond-summary}, their corresponding discrete-time Fliess operators always converge locally. Therefore, in the following sections this case will be considered in further detail.

\begin{table}[t]
\caption{Summary of convergence conditions for $F_c$ and $\hat{F}_c$}
\label{tbl:convergence-cond-summary}
\begin{center}
\vspace*{0.1in}
\begin{tabular}{lcc} \hline
 Growth Rate & $F_c$ & $\hat{F}_c$ \\ \hline
$\abs{(c,\eta)}\le K_c M_c^{|\eta|}|\eta|!$ & LC & divergent \\ \hline
$\abs{(c,\eta)}\le K_c M_c^{|\eta|}$ & GC & LC \\ \hline
$\abs{(c,\eta)}\le K_c M_c^{|\eta|}\frac{1}{|\eta|!}$ & GC (at least) & GC
\\[0.03in] \hline
\end{tabular}
\end{center}
\end{table}


\subsection{Rational discrete-time Fliess operators}
\label{sec:rational-operators}

An example of a system that can be described by a discrete-time Fliess
operator is
\allowdisplaybreaks
\begin{align} \label{eq:diff_eq_2}
\begin{split}
\left[ \begin{array}{c} z_1(N+1) \\ z_2(N+1) \end{array}\right]  = & {\;}
\left[\begin{array}{cc} 1 & 0 \\ 0 & 1  \end{array}\right] \left[
\begin{array}{c} z_1(N) \\ z_2(N) \end{array}\right] + \left[ \begin{array}{c}
1 \\ 0 \end{array}\right] \hat{u}_1(N+1) + \left[ \begin{array}{c}
0 \\ {z_1(N)} \end{array}\right] \hat{u}_2(N+1) \\
& {\;} + \left[ \begin{array}{c}
0 \\ 1 \end{array}\right] \hat{u}_1(N+1)\hat{u}_2(N+1),\\
y(N) = & {\;}
\left[\begin{array}{cc} 0 & 1 \end{array}\right] \left[
\begin{array}{c} z_1(N) \\ z_2(N) \end{array}\right],
\end{split}
\end{align}
where $\hat{u}_1$ and $\hat{u}_2$ are suitable input sequences. Letting $z_1(0)=z_2(0)=0$. Observe that
\begdi
	z_1(N+1)=z_1(N)+\hat{u}_1(N+1)
\enddi
implies that $z_1(N)=\sum_{k=1}^N \hat{u}_1(k)$. Thus, it follows that
\begin{align} \label{eq:diff_eq_3}
	\nonumber z_2(N+1)
	&= z_2(N)+z_1(N)\hat{u}_{2}(N+1)+\hat{u}_{2}(N+1)\hat{u}_{1}(N+1)\\
	& = z_2(N)+\hat{u}_{2}(N+1)z_1(N+1) \\
\nonumber & = \sum_{k_2=1}^{N+1} \hat{u}_2(k_2)\sum_{k_1=1}^{k_2}
\hat{u}_1(k_1).
\end{align}
The corresponding output is then
\begin{align*}
	y(N) = \sum_{k_2=1}^{N+1} \hat{u}_2(k_2)\sum_{k_1=1}^{k_2}\hat{u}_1(k_1) = S_{x_2x_1}[\hat{u}](N),
\end{align*}
which has the form of \rref{eq:DT-Fliess-operator-defined}. System \rref{eq:diff_eq_2} falls into the category of {\em polynomial input and
state affine} systems \cite{Sontag_79}. A simple discretization procedure can also yield discrete-time systems that are rational
functions of the inputs. Consider, for instance, the following continuous-time system 	
\begeq
\label{eq:CT-rational-Ferfera-system-realization}
	\dot{z}(t)= z(t)u(t), \;\;z(0)=0.
\endeq
For small $\Delta>0$, an Euler type approximation gives
\begin{align*}
	\tilde{z}((N+1)\Delta)&=\tilde{z}(N\Delta)+\int_{N\Delta}^{(N+1)\Delta}\tilde{z}(t)u(t)\,dt \\
	&\approx \tilde{z}(N\Delta)+\int_{N\Delta}^{(N+1)\Delta}\hspace*{-0.1in}u(t)\,dt\,\, \tilde{z}((N+1)\Delta)\\
	&=\tilde{z}(N\Delta)+\hat{u}(N+1)\,\tilde{z}((N+1)\Delta),
\end{align*}
and therefore, letting $\hat{z}(N)=\tilde{z}(N\Delta)$, observe that
\begin{equation} \label{eq:example_rational_discretization}
	\hat{z}(N+1)=(1-\hat{u}(N+1))^{-1}\hat{z}(N)
\end{equation}
In this case, $(1-\hat{u}(N+1))^{-1}$ is a rational function and fall into the following class of systems.

\begde \label{def:state-affine-rational-system} {\cite{Gray-etal_NM}}
A discrete-time state space realization is \bfem{rational input} and
\bfem{state affine} if its transition map has the form
\begin{displaymath}
	\hat{z}_i(N+1)=\sum_{j=1}^n r_{ij}(\hat{u}(N+1))\hat{z}_j(N)+s_i(\hat{u}(N+1)),
\end{displaymath}
$i=1,2,\ldots,n$, where $\hat{z}(N)\in\mathbb{R}^n$, $\hat{u}=[\hat{u}_0,\hat{u}_1,\ldots,\hat{u}_m]^T$, $r_{ij}$ and $s_i$ are rational functions, and the output map $h:\hat{z}\mapsto \hat{y}$ is linear.
\endde

The general situation is described by the following realization theorem.

\begth \label{th:rational-state-affine-realization} \cite{Gray-etal_NM}
Let $c\in\mathbb{R}\langle\langle X \rangle\rangle$ be a rational series over $X=\{x_0,x_1,\ldots,x_m\}$ with linear
representation $(\mu,\gamma,\lambda)$.
Then $\hat{y}=\hat{F}_c[\hat{u}]$ has a finite dimensional rational input and state affine realization on $B_\infty^{m+1}[0,N_f](\hat{R})$ for
any $N_f> 0$ provided $\hat{R}<\left(\sum_{j=0}^m \norm{\mu(x_j)}\right)^{-1}$, where $\norm{\cdot}$ is any matrix norm.
\endth


\subsection{Parallel interconnections and the quasi-shuffle algebra}
\label{subsec:parallel-interconnections-quasi-shuffle}

Given two continuous-time Fliess operators $F_c$ and $F_d$ with $c,d\in\allseriesLC$, the parallel interconnections as shown in Figure~\ref{fig:parallel} satisfy $F_c+F_d=F_{c+d}$ and $F_cF_d=F_{c\shuffle d}$ \cite{Fliess_81}.
In the discrete-time case, the parallel sum interconnection is characterized trivially by the addition of generating series, i.e., $\hat{F}_c+\hat{F}_d=\hat{F}_{c+d}$ due to the vector space nature of $\allseriesLC$. But the parallel product connection in this case is characterized by the so-called \emph{quasi-shuffle} product introduced in Example~\ref{exam:shuffle}. The main objective of this section is to give a description of the
quasi-shuffle algebra $H_{qsh} = (\allpoly, \circledast)$ in the context of discrete-time Fliess operators and show that rationality is preserved under the quasi-shuffle product.


\subsubsection{Quasi-shuffle algebra}
\label{ssect:QshuAlg}

The shuffle product \rref{shuffleproduct} describes the product of iterated integrals. However, it cannot account for products of iterated sums. For instance, observe that the product
\begin{equation}  \label{eq:product_summations_ex}
	\sum^{N}_{i=1}\hat{u}_1(i)\sum^{N}_{j=1}\hat{u}_2(j)
	=  \sum^N_{i = 1}\sum^i_{j = 1} \hat{u}_1(i)\hat{u}_2(j)
		+  \sum^N_{i = 1}\sum^i_{j = 1}\hat{u}_1(j)\hat{u}_2(i)
		- \sum_{i= 1}^N \hat{u}_1(i)\hat{u}_2(i),
\end{equation}
where
$\hat{u} \in B_\infty^{m+1}[0,N_f]({R})$ for suitable $R$ and $N_f$. If  $X=\{x_0,x_1,x_2\}$, then using
\rref{eq:interated-sum-approx-for-iterated-integral} it follows that \rref{eq:product_summations_ex} can be written as
\begin{align*}
	S_{x_1}[\hat{u}](N)S_{x_2}[\hat{u}](N)
	=  S_{x_1x_2}[\hat{u}](N) + S_{x_2 x_1}[\hat{u}](N)
		- \sum_{i= 1}^N \hat{u}_1(i)\hat{u}_2(i).
\end{align*}
Note that the last term $\sum_{i=1}^N\hat{u}_1(i)\hat{u}_2(i)$ does not correspond to a letter in $X$ nor to a word in $X^\ast$. Therefore, the alphabet $X$ needs to be augmented to account for this fact. Associating the input $\hat{u}_1 \hat{u}_2$ with the new letter $x_{1,2}$, one can now write
\begin{align*}
	S_{x_1}[\hat{u}](N)S_{x_2}[\hat{u}](N)
	=  S_{x_1x_2}[\hat{u}](N) + S_{x_2 x_1}[\hat{u}](N)
		+ S_{x_{1,2}}[\hat{u}](N).
\end{align*}
{{
Therefore, the general setting in which products of iterated sums are considered requires a countable alphabet. The extra letters, in addition to those in $X=\{x_0,x_1, \ldots, x_m\}$, account for all possible finite products of inputs. Recall item \ref{quasishuf} in Example \ref{exam:shuffle}, where the quasi-shuffle Hopf algebra is defined. Here the alphabet $X$ is extended to a graded commutative semigroup by defining the commutative bracket operation of letters in $X$ to be $[x_i x_j] = x_{i,j}$, which is assumed to be associative, i.e., $[[x_i x_j]x_l] = [x_i [x_jx_l]]$ for letters $x_i,x_j,x_l \in X$. Iterated brackets may therefore be denoted by $x_{i_1,\ldots,i_n}:=[[[x_{i_1}x_{i_2}]\cdots ]x_{i_n}]$. The augmented alphabet $\bar{X}$ contains $X$ as well as all finitely iterated brackets $x_{i_1,\ldots,i_n}$. The monoid of words with letters from $\bar{X}$ is denoted $\bar{X}^*$. The definition \rref{eq:interated-sum-approx-for-iterated-integral} of iterated sums has to be extended to include the additional words in $\bar{X}^*$, for instance,}}
$$
	S_{x_kx_{i_1,i_2,\ldots,i_n}}[\hat{u}](N):=\sum^N_{i = 1}\hat{u}_{k}(i)\sum^{i}_{j=1}\hat{u}_{i_1}(j)\hat{u}_{i_2}(j)\cdots \hat{u}_{i_n}(j).
$$
It follows now that the product $S_{x_1}[\hat{u}](N)S_{x_2}[\hat{u}](N)$ is encoded symbolically in terms of a quasi-shuffle product on $\bar{X}^*$
\begin{align} \label{eq:Qshuffle_2_words}
	x_1 \circledast x_2 = x_1 x_2 + x_2 x_1 - x_{1,2} \in \allpolybarX.
\end{align}
The foundation of discrete-time Fliess operator theory is the summation operator, which is used inductively in the construction of the iterated sums in \rref{eq:interated-sum-approx-for-iterated-integral}. In general, the summation operator $Z$ is defined as
\begin{equation}
\label{summation}
	Z(f)(x) := \sum_{k= 1}^{[x/\theta]} \theta f(\theta k)
\end{equation}
for a suitable class of functions $f$. It is known to satisfy the so-called Rota--Baxter relation of weight $\theta$
\cite{EFPatrasSMF_13}
\begin{equation}
\label{eq:summation_operator_product}
	Z(f)(x) Z(g)(x) = Z\bigl(Z(f)g + fZ(g)  - \theta fg\bigr)(x).
\end{equation}
This relation generalizes the integration by parts rule for indefinite Riemann integrals and provides the corresponding formula for iterated sums. Specifically, \rref{eq:product_summations_ex} corresponds to  \rref{eq:summation_operator_product} where $\theta =1$, $f=\hat{u}_1$ and $g=\hat{u}_2$.
{The quasi-shuffle product, introduced in item \ref{quasishuf} of Example \ref{exam:shuffle}, defined on $\bar{X}^*$ provides an extension of \rref{eq:Qshuffle_2_words} and \rref{eq:summation_operator_product}.} For words $\eta = \eta_{1}\cdots \eta_{{n}}$ and $\xi = \xi_{1} \cdots \xi_{{m}}$, where $\eta_i,\xi_j\in\bar{X}$, the recursive definition of the quasi-shuffle product on $\bar{X}^*$ is given by
\begeq
\label{qshu-prod-rec1}
	\eta \circledast \xi
         = \eta_{1}(\eta_{1}^{-1}(\eta) \circledast \xi) + \xi_{{1}}(\eta \circledast \xi_{1}^{-1}(\xi))
         - [\eta_{1}\xi_{1}]\big(\eta_{1}^{-1}(\eta) \circledast \xi_{1}^{-1}(\xi) \big)
\endeq
with $\emptyset \circledast \eta = \eta \circledast \emptyset = \eta$ for $\eta \in \bar{X}^*$, {and $\eta_{1}^{-1}(\cdot)$ is the left-shift operator defined in~\rref{eq:left_shift}. This implies that
\begeq \label{eq:Summation_products}
S_{\eta}[\hat{u}](N)\cdot S_{\xi}[\hat{u}](N) = S_{\eta\circledast \xi}[\hat{u}](N)
\endeq
with $\eta\circledast \xi \in \allpolybarX$. Observe that since $\abs{\eta},\abs{\xi} <\infty$, then $\supp\{\eta\circledast \xi\}$ is generated by a finite subset of $\bar{X}$.} The quasi-shuffle product $\circledast$ is linearly extended to series $c,d\in \mathbb{R}\langle\langle \bar{X} \rangle \rangle$ so that
\begdi
	c \circledast d=
	\sum\limits_{\eta,\xi \in \bar{X}^*}(c,\eta)(d,\xi)\eta \circledast \xi
	=\sum_{\nu\in \bar{X}^*} \underbrace{\sum\limits_{\eta, \xi \in \bar{X}^*}(c,\eta)(d,\xi)
	(\eta \circledast \xi,\nu)}_{\displaystyle (c\circledast d, \nu)}\nu.
\enddi
{Note that the coefficient $(\eta{\circledast}\xi,\nu)\neq 0$ only when $\nu \in X^\ast$ is such that}  $\abs{\eta}+\abs{\xi} -\min(\abs{\eta},\abs{\xi} )\le \abs{\nu} \le \abs{\eta} + \abs{\xi}$. Therefore, $(c \circledast d, \nu)$ is finite since the set $I_{\circledast}(\nu)\triangleq\{(\eta,\xi)\in \bar{X}^* \times \bar{X}^* : (\eta \circledast \xi,\nu)\not = 0\}$ is finite. {Hence, the summation defining $c\circledast d$ is locally finite, and therefore summable.} It can be shown that the quasi-shuffle product is commutative, associative and distributes over addition   {\cite{Ebrahimi-Fard-Guo_2006,Hoffman_00}}. Thus, the vector space $\mathbb{R}\langle\langle \bar{X} \rangle \rangle$ endowed with the quasi-shuffle product forms a commutative $\re$-algebra, the so-called \emph{quasi-shuffle algebra} with multiplicative identity element~$\mathbf{1}$.


\subsubsection{Rationality of the quasi-shuffle product}
\label{ssect:rationQuasiShuffle}

In this section the question of whether the quasi-shuffle product of two
rational series is again rational is addressed. In light of Definition~\ref{def:rationality}
{and the remark thereafter, it is clear that a rational series $c$ over $\bar{X}$
is also rational over a finite sub-alphabet $X_c \subset \bar{X}$. In which case,
Theorems~\ref{th:Schutzenberger} and \ref{th:stable-rational-subspace} still apply in the present setting.}
Also note that in the context of the parallel product connection
the underlying alphabets for the generating series of $\hat{F}_c$ and $\hat{F}_d$ are always the same since the inputs are identical. But there is no additional complexity introduced if the alphabets are allowed to be distinct. So let $X_c, X_d \subset \bar{X}$ be finite sub-alphabets of {$\bar{X}$ corresponding to the generating series $c$ and $d$ and with cardinalities $N_c$ and $N_d$, respectively}. Define
$[X_c X_d]= \{[x^c_i x^d_j]: x^c_i\in
X_c,x^d_j\in X_d, i=1,\ldots,N_c, j=1\ldots,N_d\}$.
The main theorem of the section is given first.

\begth \label{th:Qshuffle_rationality}
{Let $c,d \in \allseriesbarX$ be rational series with underlying finite alphabets $X_c,X_d\subset \bar{X}$, then $e=c \circledast d \in \allseriesbarX$ is rational with underlying alphabet $X_e= X_c\cup X_d \cup [X_c X_d]\subset \bar{X}$.}
\endth

\begin{proof}
In light of \rref{qshu-prod-rec1}, the series $e=c\circledast d$ is clearly defined over the finite alphabet $X_e$.
{Therefore, a stable finite dimensional vector space $V_e$ is constructed which contains $e$ in
order to apply Theorem~\ref{th:stable-rational-subspace}.}
Since $c$ and $d$ are both rational,
let $V_c$ and $V_d$ be stable finite
dimensional vector subspaces of $\allseriesXc$ and $\allseriesXd$
containing $c$ and $d$,
respectively. Let $\{\bar{c}_i\}_{i=1}^{n_c}$ and $\{\bar{d}_j\}_{j=1}^{n_d}$
denote their corresponding bases.
Define
\begdi
V_{e}={\rm span}_{\re}\{\bar{c}_i\circledast\bar{d}_j\;:\; i=1,\hdots,n_c,\;\; j=1,\hdots,n_d  \}.
\enddi
Clearly, $V_{e}\subset\allseriesXe$ is finite dimensional. If one
writes
\begdi
c=\sum^{n_c}_{i=1}\alpha_i \bar{c}_i,  \;\; d=\sum^{n_d}_{j=1} \beta_j \bar{d}_j,
\enddi
it then follows directly that
\begdi
e=c\circledast d = \sum^{n_c,n_d}_{i,j=1}\alpha_i\beta_j\;
\bar{c}_i\circledast\bar{d}_j\in V_e.
\enddi
So it only remains to be shown that $V_{e}$ is stable. Observe from
\rref{qshu-prod-rec1} that
for any $x\in X_e$ the left-shift operator acts on the quasi-shuffle product as
\begin{equation} \label{eq:Qshuffle_on_shift_op}
x^{-1}(\eta \circledast \xi) = x^{-1}(\eta) \circledast \xi
+ \eta \circledast x^{-1}(\xi) +
\delta_{x,[x_i x_j]}( x_i^{-1}(\eta)\circledast x_j^{-1}(\xi)),
\end{equation}
where $\eta=x_i\eta',\xi = x_j \xi' \in \bar{X}^\ast$ and $\delta_{x,[x_i x_j]}
=1 $
if
$x=[x_ix_j]$ and $0$ otherwise. Writing $c= (c,\emptyset) + \sum_{i=0}^{N_c}
x^c_i
\; (x^c_i)^{-1}(c)$ and $d= (d,\emptyset) + \sum_{i=0}^{N_d} x^d_i \;
(x^d_i)^{-1}(d)$
and using the bilinearity of the quasi-shuffle, it follows that
\begdi
x^{-1}(e) = {x^{-1}(c) \circledast d
+  c \circledast x^{-1}(d)}  +\sum_{i,j=0}^{N_c,N_d}
\delta_{x,[x^c_i x^d_j]}( (x^c_i)^{-1}(c)\circledast
(x^d_j)^{-1}(d)).
\enddi
But since $V_c,V_d \subset \allseriesXe$ are stable vector spaces by assumption, it is immediate that {$(x^c_i)^{-1}(c) \in V_c$ and $(x^d_j)^{-1}(d) \in V_d$, and therefore $x^{-1}(e)\in V_e$ as well. It then follows that $V_{e}$ is a stable vector space, and hence $e$ is rational.}
\end{proof}

The following corollary describes the generating series for the parallel product connection
in the context of rational series.

{\begco \label{co:Qshuffle_rationality}
If $c,d \in \allseriesbarX$ are rational series with underlying finite alphabets $X_c,X_d\subset \bar{X}$, then $\hat{F}_c \hat{F}_d = \hat{F}_{c\circledast d}$ with $e=c \circledast d
\in \allseriesXeRat$, where $X_e= X_c\cup X_d \cup [X_c X_d]$.
\endco
\begin{proof} From \rref{eq:Summation_products} the product connection of two operators as in \rref{eq:DT-Fliess-operator-defined} is
\begin{align*}
\hat{F}_c[\hat{u}](N) \hat{F}_d[\hat{u}](N) &= \sum\limits_{\eta\in X_c^*}(c,\eta)S_{\eta}[\hat{u}](N)\cdot\sum\limits_{\xi\in X_d^*}(d,\xi)S_{\xi}[\hat{u}](N)\\
&=  \sum\limits_{\eta\in X_c^*,\xi\in X_d^*}(c,\eta)(d,\xi)S_{\eta\circledast\xi}[\hat{u}](N)\\
&=  F_{c\circledast d}[\hat{u}](N) =: F_e[\hat{u}](N).
\end{align*}
Here $e\in\allseriesXeRat$ since by Theorem \ref{th:Qshuffle_rationality} the quasi-shuffle preserves rationality.
\end{proof}}

The following lemma will be used in the final example of this section.

\begle \label{le:Qshuffle_identity}
For any $i,j\ge 0$
\begin{equation} \label{eq:Qshuffle_identity_x1}
	x_1^i \circledast x^j_1 = \sum_{k=0}^{\min\{i,j\}}
	\binom{i+j-2k}{\min\{i,j\}-k}
	x^k_{1,1} \shuffle x_1^{i+j-2k}.
\end{equation}
\endle

\begin{proof}
Without loss of generality assume $i\ge j$. The identity is proved by
induction over $i+j$. The cases for $i+j=0,1$ are trivial. Assume
\rref{eq:Qshuffle_identity_x1} holds up to some fixed $i+j$. Using
\rref{qshu-prod-rec1} compute
\begin{align*}
	x_1^i \circledast x^{j+1}_1
	= x_1\left( x_1^{i-1} \circledast x^{j+1}_1\right) + x_1\left(
	x_1^{i} \circledast x^{j}_1\right) + x_{1,1}\left( x_1^{i-1} \circledast x^{j}_1\right).
\end{align*}
By the induction hypothesis and since $i\le j$,
\begin{align*}
	\lefteqn{x_1^i \circledast x^{j+1}_1
	= {\;}  \sum_{k=0}^{i-1} \binom{i+j-2k}{i-1-k} x_1(x^k_{1,1} \shuffle x_1^{i+j-2k}) } \\
	&+ \sum_{k=0}^{i} \binom{i+j-2k}{i-k} x_1 (x^k_{1,1} \shuffle x_1^{i+j-2k})
	 + \sum_{k=0}^{i-1} \binom{i+j-1-2k}{i-1-k} x_{1,1} (x^k_{1,1} \shuffle x_1^{i+j-1-2k} )\\
	= &{\;}  \sum_{k=0}^{i-1} \binom{i+j-2k}{i-1-k} x_1 (x^k_{1,1} \shuffle x_1^{i+j-2k}) \\
	& + \sum_{k=0}^{i-1} \binom{i+j-2k}{i-k} x_1 (x^k_{1,1} \shuffle x_1^{i+j-2k})
	+ \binom{j-i}{0} x_1(x_{1,1}^k\shuffle x^{j-i})\\
	& + \sum_{k=1}^{i-1} \binom{i+j-2k+1}{i-k} x_{1,1}(x^k_{1,1} \shuffle x_1^{i+j-2k+1})
	+ \binom{j-i+1}{0} x_1(x_{1,1}^{k-1}\shuffle x^{j-i+1})\\
	= &{\;} (x^k_{1,1} \shuffle x_1^{i+j-2k}) + \binom{i+j+1}{i}x^{i+j+1}_1
	+ \sum_{k=0}^{i-1} \binom{i+j-2k+1}{i-k} x_1 (x^k_{1,1} \shuffle x_1^{i+j-2k}) \\
	& + \sum_{k=1}^{i-1} \binom{i+j-2k+1}{i-k} x_{1,1}(x^k_{1,1} \shuffle x_1^{i+j-2k+1})\\
	=& {\,} \sum_{k=0}^{i} \binom{i+j-2k+1}{i-k} x^k_{1,1} \shuffle x_1^{i+j-2k+1}.
\end{align*}
This complete the proof since it was assumed that $\min\{i,j\} = i$.
\end{proof}

\begex
Let $X=\{x_1\}$ and
consider the rational series $c=x_1^\ast:=\sum_{k\geq 0} x_1^k$. It can be shown directly that
\begin{align} \label{eq:ex_CT_product_x1ast}
x_1^{\ast} \shuffle x_1^{\ast} = \sum_{n=0}^\infty \sum_{i=1}^n
 \binom{n}{i} x_1^{n} = \sum_{\eta\in X^\ast} 2^{\abs{\eta}} \eta,
\end{align}
using the identity $x_1^i\shuffle x_1^j = \binom{i+j}{i} x_1^{i+j}$ \cite{Wang_90}. Since the
shuffle product is known to preserve rationality, it follows from
Theorem \ref{th:Schutzenberger} that $x_1^{\ast} \shuffle x_1^{\ast}$ must have a linear
representation $(\mu,\gamma,\lambda)$, in this case $\mu(\eta)=
2^{\abs{\eta}}$ and $\gamma=\lambda=1$. This is easily verified by setting $z_i=F_c[u]$,
which gives the bilinear state space realization $\dot{z}_i=z_iu$, $y_i=z_i$. Then the
parallel product connection $y=y_1y_2=F^2_c[u]=z$ has the realization
$\dot{z} = 2 z u$, $y=z$.
One can confirm using iterated Lie derivatives that the
generating series for this system is exactly \rref{eq:ex_CT_product_x1ast}.
\endex

\begex
The goal now is to produce the quasi-shuffle analogue of \rref{eq:ex_CT_product_x1ast}.
Note here that $X=\{x_1,x_{1,1}\}$.
Using Lemma~\ref{le:Qshuffle_identity} it follows that
\begin{align*}
	x_1^{\ast} \circledast x_1^{\ast} =& {\,} \sum_{i,j=0}^\infty x_1^i \circledast x_1^j
	=  \sum_{n=0} \sum_{i+j=n} \sum_{k=0}^{\min\{i,j\}} \binom{i+j-2k}{\min\{i,j\}-k} x_{1,1}^k \shuffle x_1^{i+j-2k}.
\end{align*}
For fixed $k',n' \in \nat \cup \{0\}$,  let $\eta\in X^{n'}$ be such that $\abs{\eta}_{x_{1,1}} = k'$, where $|\eta|_{x_{i}}$ denotes the number of times the letter  $x_i\in X$ appears in $\eta\in X^\ast$. It follows that 
the coefficient $(x_1^{\ast} \circledast x_1^{\ast},\eta)$ is written as
\begin{align} \label{eq:coef_Qshuffle1}
(x_1^{\ast} \circledast x_1^{\ast},\eta)  
= & {\,}  \sum_{n=0} \sum_{i+j=n} \sum_{k=0}^{\min\{i,j\}} \binom{i+j-2k}{\min\{i,j\}-k} (  x_{1,1}^k \shuffle x_1^{i+j-2k} , \eta ).
\end{align}
Notice first that the number of elements in $\supp(x_{1,1}^k \shuffle x_{1}^{i+j-2k})$ is $\binom{i+j-k}{k}$, and secondly that it implies that $(x_{1,1}^k \shuffle x_{1}^{i+j-2k},\eta) = 1$ when $k=k'$ and $i+j=n'+k'$, otherwise $(x_{1,1}^k \shuffle x_{1}^{i+j-2k},\eta) = 0$. The former statement is related to the fact that $ Q^N = \sum_{i_1+\cdots +i_m = N} q_1^{i_1} \shuffle \cdots \shuffle q_m^{i_m} $ for an arbitrary alphabet $ Q=\{q_1,\ldots,q_m\}$ \cite{Duffaut_Espinosa_09}.  The coefficient $(x_1^{\ast} \circledast x_1^{\ast},\eta) $ can be further developed as
\begin{align*}
(x_1^{\ast} \circledast x_1^{\ast},\eta)
= &{\,} \sum_{i+j=n'+k'}\binom{n'-k'}{\min\{i,j\}-k'}
= \sum_{i+j=n'+k'} \binom{n'-k'}{\min\{i-k',j-k'\}} \\
= &{\,} \sum_{i=k'}^{n'} \binom{n'-k'}{ \min\{i-k',n'-i\}}
= \sum_{i=k'}^{n'} \binom{n'-k'}{i-k'} =\sum_{i=0}^{n'-k'} \binom{n'-k'}{i}  = 2^{\abs{\eta}_{x_1}}.
\end{align*}
Thus, one can write
\begin{align} \label{eq:Qshuffle_Ferfera_with_coefficients}
	x_1^{\ast} \circledast x_1^{\ast}
	= &{\,}\sum_{\eta\in X^\ast} (x_1^{\ast} \circledast x_1^{\ast},\eta) \eta =  \sum_{\eta\in X^\ast} 2^{\abs{\eta}_{x_1}} \eta.
\end{align}
In light of Theorem~\ref{th:Qshuffle_rationality}, the series $x_1^{\ast} \circledast x_1^{\ast}$ must be rational.
In particular, a straightforward linear representation for $x_1^{\ast} \circledast x_1^{\ast}$ can be obtained due to \rref{eq:Qshuffle_Ferfera_with_coefficients}. That is, $(\mu ,\gamma,\lambda)$ is identified as $\mu(\eta) = 2^{\abs{\eta}_{x_1}}$ and
$\lambda=\gamma =1$. Finally, a direct computation confirms that
\begin{align*}
x_1^{\ast} \circledast x_1^{\ast} = &{\;} \emptyset + 2x_1 + x_{1,1} + 4x_1^2+2x_1x_{1,1}+2x_{1,1}x_1+x_{1,1}^2 +  8x_1^3 + 4 x_1^2 x_{1,1} \\
&{\;}  +4 x_1 x_{1,1} x_1 + 4x_{1,1} x_1^2 +2x_1 x_{1,1}^2 + 2 x_{1,1} x_1 x_{1,1} + 2 x_{1,1}^2 x_1 + x_{1,1}^3 + 16 x_1^4 \\
&{\;} + 8 x_1^3 x_{1,1} + 8 x_1^2 x_{1,1} x_1 + 8 x_1 x_{1,1} x_1^2 + 8 x_{1,1} x_1^3 + 4 x_1^2 x_{1,1}^2 +  4 x_1 x_{1,1} x_1 x_{1,1} \\
&{\;} + 4 x_1 x_{1,1}^2 x_1  + 4 x_{1,1} x_1^2 x_{1,1} +
4 x_{1,1} x_1 x_{1,1} x_1 + 4 x_{1,1}^2 x_1^2  + 2 x_1 x_{1,1}^3  + 2 x_{1,1} x_1 x_{1,1}^2 \\
&{\;}  + 2 x_{1,1}^2 x_1 x_{1,1} + 2 x_{1,1}^3 x_1 + x_{1,1}^4 +\cdots,
\end{align*}
where it is clear that \rref{eq:Qshuffle_Ferfera_with_coefficients} holds. \endex

%

\begin{acknowledgement}
The third author was supported by grant SEV-2011-0087 from the Severo Ochoa Excellence Program at the Instituto de Ciencias Matem\'{a}ticas in Madrid, Spain. This research was also supported by a grant from the BBVA Foundation.
\end{acknowledgement}



\end{document}